\documentclass[12pt]{amsart}
\usepackage{amsmath}
\usepackage{amsfonts}
\usepackage[mathscr]{euscript}
\usepackage[arrow,matrix,curve]{xy}
\usepackage{xyadds} 
\usepackage{array}
\def\mpar#1{\relax}
\input arthead

\DeclareMathAlphabet      {\mathsfi}{OT1}{cmss}{m}{sl}
\DeclareMathAlphabet      {\mathsfb}{OT1}{cmss}{bx}{n}

\makeatletter
\def\skd{\mathsfi} 
\def\sketchname{\mathsfi}
\def\catname{\mathbf}

\def\pd{.\penalty\exhyphenpenalty} 

\let\prd\prod 

\def\fake{\hspace*{.01pt}} 

\def\logopn#1{\expandafter\def\csname#1\endcsname{%
\mathop{\mathsf{#1\fake}}\nolimits}}
\def\logord#1{\expandafter\def\csname#1\endcsname{%
\mathord{\mathsf{#1\fake}}}}
\def\logvary#1{\expandafter\def\csname#1\endcsname{%
\mathop{\mathsfi{#1\fake}}\nolimits}}
\def\csspace#1{\expandafter\def\csname#1\endcsname{%
\mathop{\mathsfb{#1\fake}}\nolimits}}
\def\skname#1{\expandafter\def\csname#1\endcsname{%
\mathop{\mathsfi{#1\fake}}\nolimits}}
\def\frakop#1{\expandafter\def\csname#1\endcsname{%
\mathop{\mathfrak{#1\fake}}\nolimits}}

\csspace{BinProd}
\csspace{Cat}
\csspace{CCC}
\csspace{FinProd}
\csspace{FinLim}

\frakop{C}
\frakop{M}
\frakop{U}

\logopn{AmbSp}
\logopn{Arrows}
\logopn{Arr}
\logopn{BsDiag}
\logopn{CatTh}
\logopn{CatStruc}
\logopn{Cod}
\logopn{Cones}
\logopn{Constants}
\logopn{ConstSp}
\logopn{Ddcn}
\logopn{Diag}
\logopn{Diagrams}
\logopn{Dom}
\logopn{Expr}
\logopn{Fillin}
\logopn{FLTh}
\logopn{FLUnivMod}
\logopn{Form}
\logopn{FPTh}
\logopn{Func}
\logopn{Graph}
\logopn{Id}
\logopn{Incl}
\logopn{Init}
\logopn{InputTypes}
\logopn{Inp}
\logopn{Length}
\logopn{Left}
\logopn{LimCone}
\logopn{Lim}
\logopn{Lin}
\logopn{LinTh}
\logopn{LinUnivMod}
\logopn{List}
\logopn{Nodes}
\logopn{Oprns}
\logopn{Outp}
\logopn{Pbl}
\logopn{Pbr}
\logopn{PEqAr}
\logopn{PPEqAr}
\logopn{Proj}
\logopn{Right}
\logopn{Rng}
\logopn{ShpGr}
\logopn{Sk}
\logopn{SynCat}
\logopn{Th}
\logopn{TypeList}
\logopn{TypeSet}
\logopn{Types}
\logopn{Type}
\logopn{UndCat}
\logopn{UndSk}
\logopn{UndGr}
\logopn{UnivMod}
\logopn{VarList}
\logopn{VarSet}
\logopn{Var}
\logopn{Vbl}
\logopn{Vertex}

\logvary{claim}
\logvary{claimcon}
\logvary{f}
\logvary{hyp}
\logvary{hypcon}
\logvary{verif}
\logvary{wksp}

\skname{E}
\skname{F}
\skname{S}
\skname{T}

\def\arr{\mathord{{\sf ar\fake}}}
\logord{arrow}
\logord{base}
\logord{bot}
\logord{comp}
\logord{cone}
\logord{csource}
\logord{ctarget}
\logord{curry}
\logord{ebot}
\logord{econe}
\logord{efid}
\logord{eq}
\logord{equ}
\logord{esoco}
\logord{etaco}
\logord{etop}
\logord{eufid}
\logord{efia}
\logord{ev}
\logord{exp}
\logord{inc}
\logord{fia}
\logord{fid}
\logord{fs}
\logord{fss}
\logord{id}
\logord{lam}
\logord{lass}
\logord{lfac}
\logord{lproj}
\logord{lunit}
\logord{mfac}
\logord{modpon}
\logord{n}
\logord{ob}
\logord{ppair}
\logord{prod}
\logord{proj}
\logord{rass}
\logord{rfac}
\logord{runit}
\logord{rproj}
\logord{s}
\logord{soco}
\logord{soob}
\logord{source}
\logord{target}
\logord{ta}
\logord{ter}
\logord{taco}
\logord{top}
\logord{tsource}
\logord{ttarget}
\logord{twovf}
\logord{unit}
\logord{ufid}
\logord{v}
\logord{vertex}

\def\Fa{\cor\F}

\def\specialsubsection{\subsubsection}

\def\@listI{\leftmargin\leftmargini \parsep 5pt plus 2.5pt minus 1pt\topsep
10pt plus 4pt minus 6pt\itemsep 3pt plus 2pt minus 1pt}
\let\@listi\@listI
\@listi

\def\Biggg{\bBigg@{5}}
\def\Bigggg{\bBigg@{7}}
\def\biggggg{\bBigg@{8}}
\def\Biggggg{\bBigg@{9}}

\def\brspace{\hspace{3em}}

\makeatother

\renewenvironment{squishlist}
{\begin{list}{\arabic{lister})}%
{\usecounter{lister}\parsep0pt\itemsep 0pt plus 2pt
\listparindent 1.5em}}
{\end{list}}
\renewcommand{\sqlist}{\begin{squishlist}}
\renewcommand{\esqlist}{\end{squishlist}\noindent}

\mathchardef\eqs="603D
\mathchardef\lra="622C
\mathchardef\xord="0202 

\def\Frac#1#2{{\def\arraystretch{1.2}
\begin{array}{c}{#1}\\[6pt] \hline
\\[-12pt]{#2}\end{array}}}

\def\Rule#1#2#3#4{\mbox{\rm #3}%
\hspace*{1.5em}{\Frac{\kern-2pt #1}%
{\kern2pt #2}}
\hspace*{1.5em} \parbox{1.8in}{#4}}

\def\rtwo#1#2{\ar[r]<1ex>^<>(.5){#1}\ar[r]<-1ex>_<>(.5){#2}}

\def\epf{$\Box$\par\addvspace{\medskipamount}}
\def\lst#1#2#3{#1_{#2},\ldots,#1_{#3}}

\def\ceq{\mathrel{\colon\hspace{-2pt}\eqs}}
\def\ceqv{\mathrel{\colon\hspace{-2pt}\lra}}

\def\astep{\subsubsection{Translation as a construction}}
\def\bstep{\subsubsection{Expression as actual factorization}}

\def\crs#1#2{\prd_{i=1}^{#2}{#1}_i}
\def\typevar{\gamma}
\def\ab#1{\Length\left[#1\right]}

\def\wider{\spreaddiagramcolumns{.5pc}}

\def\taller{\spreaddiagramrows{.5pc}}

\def\pf{\specialsubsection*{Proof}}
\makeatletter
\def\mld#1$${\null\,\vcenter\bgroup\def\\{\cr&}\openup9pt\m@th
\ialign\bgroup\strut\hfil$\displaystyle{##}$%
&$\displaystyle{{}##}$\hfil\crcr #1
\crcr\egroup\egroup\,$$}
\makeatother

\def\mk{\underline}
\def\ds#1#2#3{#1^{#2}_{#3}}
\def\dsu#1#2#3{\mk{#1}^{#2}_{#3}}

\def\cor#1{\raisebox{1ex}{$\scriptstyle\ulcorner$}%
{\hspace{-0.2em}#1}\raisebox{1ex}{$\scriptstyle\urcorner$}}

\def\scr{\mathscr}
\def\Psc{\mathord{\scr{P}}}

\def\Ssc{\mathord{\scr{S}}}

\def\Usc{\mathord{\scr{U}}}

\makeatletter
\def\biggg{\bBigg@4}
\def\Biggg{\bBigg@{4.5}}
\def\bigggg{\bBigg@{5.0}}
\makeatother

\def\mpar#1{\relax}
\def\pf{\begin{proof}}
\def\epf{\end{proof}}
\begin{document}

\title[Graph-based logic and sketches II]{Graph-based logic and
sketches II:\\ Finite-product categories and equational logic\\
(Preliminary Report)}

\author{Atish Bagchi}

\address{\hskip-\parindent Atish Bagchi\\
226 West Rittenhouse Square\\
\#702\\
Philadelphia, PA 19103}
\email{atish@math.upenn.edu}
\thanks{Research at MSI is supported in part by NSF
grant DMS-9022140.}

\author{Charles Wells}

\address{\hskip-\parindent Charles Wells\\
Department of Mathematics\\
Case Western Reserve University \\ University Circle\\
Cleveland, OH 44106-7058,  USA}
\email{cfw2@po.cwru.edu}

\begin{abstract}
It is shown that the proof theory for sketches and forms
provided in \cite{logstr} is strong enough to produce all the
theorems of the entailment system for multisorted equational
logic provided in \cite{gogmes}.\end{abstract}

\maketitle


\section{Introduction}\label{mel} In
\cite{gensk} the second author introduced the notion of
\textbf{form}, a graph-based method of specification of
mathematical structures that generalizes Ehresmann's
sketches. In \cite{logstr}, the authors produced  a
structure for forms which provides a uniform proof theory based
on finite-limit constructions for many types of forms, including
all types of sketches and also forms that can specify
higher-order structures in cartesian closed categories and
toposes (among many others).  The parameter in the proof theory
that determines the types of constructions that can be made is
the \textbf{constructor space}.  For example, the constructor
space for cartesian closed categories (with specified structure)
is the finite-limit theory $\CCC$ for cartesian closed
categories.  In particular for the concerns of the present
paper, the constructor space for structures that can be specified
by finite products is  a finite-limit theory $\FinProd$ for
categories with specified finite products.  This theory is
described explicitly in \cite{logstr}.

Each finite product form $\F$ is
given by a \textbf{syntactic category} denoted by
$\SynCat[\FinProd,\F]$. The logical structure in \cite{logstr}
identifies a statement as a \textbf{potential factorization} in
$\SynCat[\FinProd,\F]$, which is a diagram of the form
\begin{equation}\label{potf}\xymatrix{ & \hyp
\ar[d]^{\claimcon}\\ \claim \ar[r]_{\hypcon} &
\wksp}\end{equation} and the theorem that the given statement is
true as an \textbf{actual factorization}
\begin{equation*}\xymatrix{ & \hyp \ar[d]^{\claimcon}\\ \claim
\ar[ur]^{\verif} \ar[r]_{\hypcon} & \wksp}\end{equation*} of the
diagram~\eqref{potf}.

In \cite{gogmes}, Goguen and Meseguer
produced a sound and complete entailment system for multisorted
equational logic. In this paper, we verify that the theorems of
that logic for a particular signature and equations all occur as
actual factorizations in $\SynCat[\FinProd,\F]$, where $\F$ is a
$\FinProd$ form induced (in a manner to be described) by the
given signature and equations.  We also compare the expressive
powers of these two systems.

\section{Preliminaries}

\subsection{Lists} Given a set $A$, $\List[A]$ denotes the set
of lists of elements of $A$, including the empty list.  The
$k$th entry in a list $w$ of elements of $A$ is denoted by
$w_k$ and the length of $w$ is denoted by $\ab{w}$.  The \textbf{range}
of $w$, denoted by $\Rng[w]$, is the set of elements of $A$
occurring as entries in $w$.  If $f:A\to B$ is a function,
$\List[f]:\List[A]\to\List[B]$ is by definition $f$ ``mapped
over'' $\List[A]$:  If $w$ is  a list of elements of $A$, then
the $k$th entry of $\List[f](w)$ is by definition $f(w_k)$. This
makes $\List$ a functor from the category of sets to itself.

\subsection{Signatures}

\subsubsection{Expressions and terms}\label{termexp}

In the description that follows
of the terms and equations for a signature,
we use a notation that
specifies the variables of a term or equation explicitly. In
particular, one may specify
variables that do not actually appear in the expression. For
this reason, the formalism we introduce in the definitions below
distinguishes an \textit{expression\/} such as $f(x,g(y,x),z)$
from  a \textit{term\/}, which is an expression together with a
specified set of typed variables;  in this case that set could
be for example $\{x,y,z,w\}$.  This formalism is equivalent to
that of \cite{gogmes}.

\defn\label{sigf}

A pair $(\Sigma,\Omega)$ of sets together with two functions
$\Inp:\Omega\to\List[\Sigma]$ and $\Outp:\Omega\to\Sigma$
is called a \textbf{signature}.
Given a signature $\Ssc\ceq(\Sigma,\Omega)$,
elements of $\Sigma$ are called the \textbf{types} of $\Ssc$ and the
elements of $\Omega$ are called the \textbf{operations} of $\Ssc$.\edefn

\notat Given a signature $\Ssc=(\Sigma,\Omega)$,
we will denote the set $\Sigma$
of types by $\Types[\Ssc]$ and the set $\Omega$ of operations by
$\Oprns[\Ssc]$. For any $f\in\Omega$, the list $\Inp[f]$ is called
the \textbf{input type
list} of $f$ and the type $\Outp[f]$ is the \textbf{output type} of $f$.
\enotat

\rem\label{inptlrem} The input type list of $f$ is usually
called the \textbf{arity} of $f$, and the output type of $f$ is
usually called simply the \textbf{type} of $f$.\erem

\defn\label{constdef}
An operation $f$ of a signature $\Ssc$ is called a \textbf{constant}
if and only if $\Inp[f]$ is the empty list.\edefn

\defn A type $\typevar$ of a signature $\Ssc$ is said to be
\textbf{inhabited}
if and only if either \alist
\item there is a constant of output type $\typevar$ in $\Ssc$,
or
\item
there is an operation $f$
of output type $\typevar$ for
which every type in $\Inp[f]$ is inhabited.
\ealist
The type $\typevar$ is said to be \textbf{empty} if and only if
it is not inhabited.
\edefn

\subsection{Terms and equations}

In this section, we define the terms and equations
of a given signature.

\subsubsection{Assumptions}\label{assu}
In these definitions, we
make the following assumptions, useful for bookkeeping purposes.
\blist{A}\item\label{ordtyp}
We assume that we are given a signature $\Ssc$ for
which
$\Types[\Ssc]=\sb{\sigma^i}{i\in I}$ for some ordinal $I$.
\item\label{doubleind} For each $i\in I$, we assume there is an indexed set
$\Vbl[\sigma^i]\ceq \sb{x_j^i}{j\in\omega}$ whose elements are by
definition \textbf{variables of type
$\sigma^i$}. In this setting, $x_j^i$ is the $j$th variable of type
$\sigma^i$.

\item\label{orddd} The set of variables is
ordered by defining $$x_j^i\lt x_l^k \ceqv \begin{cases}\text{either} & i\lt
k\\ \text{ or }&i=k \text{ and } j\lt l\end{cases}$$ \elist

We also define $\Vbl[\Ssc]\ceq
\union_{i\in\omega}\Vbl[\sigma^i]$.

\defn For any type
$\tau$, an \textbf{expression of type} $\tau$ is defined recursively as
follows. \blist{Expr} \item A variable of type $\tau$ is an expression of type
$\tau$. \item If $f$ is an operation with $\Inp[f]=(\gamma^i\mid
i\in 1\twodots n)$ and $\Outp[f]=\tau$, and $\lstsb{e}{1}{n}$
is a list of expressions for which each $e_i$ is of type
$\gamma^i$, then $f\lstsb{e}{1}{n}$ is an expression of type
$\tau$.\elist\edefn

\notat\label{typenot}
The type of a variable $x$ is
denoted by $\Type[x]$, so that in the notation of~\ref{assu},
$\Type[x_j^i]=\sigma^i$. This notation will be extended to
include lists and sets of variables: If
$W\ceq\{x^1_2,x^1_3,x^2_1,x^3_2\}$, then
$\Type[W]=\sigma^1\xord\sigma^1\xord\sigma^2\xord\sigma^3$.
(Note that this depends on the ordering given by
A.\ref{ordtyp}.)  The type of an expression $e$ will be denoted
by $\Type[e]$. Thus the function $\Type$ is overloaded: it may
be applied to variables, sets of variables, or expressions, and
will in the following be applied to terms and equations as well.
\enotat

\defn\label{def234} For a given expression $e$,
the list of variables in $e$, in order of appearance in
$e$
from left to right, counting repetitions, is called
the \textbf{variable list} of $e$, denoted by $\VarList[e]$.
$\Rng\bigl[\VarList[e]\,\bigr]$, the set of distinct
variables occurring in $e$, is called
the \textbf{variable
set} of $e$ and denoted by $\VarSet[e]$.
The list
$(\List[\Type])\bigl[\VarList[e]\,\bigr]$
is called the \textbf{type list} of $e$, denoted by $\TypeList[e]$.
Thus if the $k$th entry of
$\VarList[e]$ is $x_j^i$, then the $k$th entry of $\TypeList[e]$
is $\sigma^i$.
The set $\Rng\bigl[\TypeList[e]\,\bigr]$, which is
the set of distinct types occurring in
$e$, is called the
\textbf{type set} of $e$, denoted by $\TypeSet[e]$.
\edefn

\exam\label{expl}
Let $e$ be the expression $f(x,g(y,x),z)$.
If $x$ and $y$ are variables of type $\typevar$
and $z$ is of type $\tau$, then the variable list of $e$
is $(x,y,x,z)$, the variable set is
$\{x,y,z\}$, the type list is $(\typevar,\typevar,\typevar,\tau)$, and
the type set is $\{\typevar,\tau\}$.  Using the notation of
A.\ref{doubleind} and supposing $\typevar=\sigma^1$,
$\tau=\sigma^2$, $x=x_1^1$, $y=x^1_2$ and $z=x^2_1$,
we have $e=f(x_1^1,g(x_2^1,x_1^1),x_1^2)$ and
the following statements hold:
\begin{align*}
\VarList[e]&=
(x_1^1,x_2^1,x_1^1,x_1^2)  \\
\VarSet[e]&=\{x_1^1,x_2^1, x_1^2\} \\
\TypeList[e]&= (\sigma^1, \sigma^1,\sigma^1,\sigma^2) \\
\TypeSet[e]&=\{\sigma^1,\sigma^2\}
\end{align*}
\eexam

\defn\label{termdeff}
A \textbf{term} $t$ for a signature $\Ssc$ is determined by the
following:
\blist{TD}\item A set
$\Var[t]$ of typed variables.  (It is a set, not a list, but it is
ordered by the ordering of
A.\ref{orddd} in~\ref{assu}.)
\item An
expression $\Expr[t]$.
\item A type $\Type[t]\in\Types[\Ssc]$.\elist
These data must satisfy the following requirements:\blist{TR}
\item $\VarSet\bigl[\Expr[t]\,\bigr]\includedin \Var[t]$.
\item $\Type[t] = \Type\bigl[\Expr[t]\,\bigr]$.\elist
\edefn
\notat A given
term $t$ will be represented as the list $$(\Expr[t],\Var[t],\Type[t])$$
\enotat

\defn\label{InputTypes}
Let $t$ be a term.
The list $\InputTypes[t]$ is defined to be the
list whose $i$th entry is the
type of the $i$th variable in
$\Var[t]$ using the ordering given by
A.\ref{orddd} in~\ref{assu}.
Thus if the $k$th entry of
$\Var[t]$ is $x_j^i$, then the $k$th entry of $\InputTypes[t]$
is $\sigma^i$.
Observe that there are no repetitions in $\Var[t]$ but
there may well be repetitions in $\InputTypes[t]$.
\edefn

\exam
Let $e=f\left(x_1^1,g(x_2^1,x_1^1),x_1^2\right)$
as in Example~\ref{expl}, and suppose $\Outp[f]=\sigma^5$.
Then there are many terms $t$ with $\Expr[t]=e$, for example
$$t_1\ceq\left(e,\{x_1^1,x_2^1,x_1^2\},\sigma^5\right)$$ and
$$t_2\ceq\left(e,\{x_1^1,x_2^1,x_3^1,x_1^2,x_5^7\},\sigma^5\right)$$
We have
$\Type[t_1] = \Type[t_2] = \sigma^5$ and (for example)
$$\InputTypes[t_1]=(\sigma^1,\sigma^1,
\sigma^2)$$ and
$$\InputTypes[t_2]=(\sigma^1,\sigma^1,
\sigma^1,\sigma^2,\sigma^7)$$\eexam

\defn An \textbf{equation} $E$ is determined by
a set $\Var[E]$ of typed variables (ordered by our convention)
and two expressions $\Left[E]$,
$\Right[E]$, for which
\blist{ER}
\item\label{firster}
$\Type\bigl[\Left[E]\,\bigr]=\Type\bigl[\Right[E]\,\bigr]$.
\item $\VarSet\bigl[\Left[E]\,\bigr]
\union\VarSet\bigl[\Right[E]\,\bigr]\includedin
\Var[E]$.\elist
\edefn
\notat
We will write $e=_Ve'$ to denote an equation $E$ with
$V=\Var[E]$, $e=\Left[E]$ and $e'=\Right[E]$.
The notation $\Type[E]$ will denote
the common type of $\Left[E]$ and $\Right[E]$.
\enotat

\exam\label{eqex2}
Let $e$ be the expression $f\left(x_1^1,g(x_2^1,x_1^1),x_1^2\right)$
of Example~\ref{expl}.  Let $e'\ceq g(x_2^1,x_3^1)$.  Then there are
many equations with $e$ and $e'$ as left and right sides, for
example:
\begin{equation}\label{mcee}
E_1 \ceq f\left(x_1^1,g(x_2^1,x_1^1),x_1^2\right)
=_{\{x_1^1,x^1_2,x_3^1,x_1^3\}}
g(x_2^1,x_3^1)
\end{equation} and
\begin{equation}
E_2 \ceq f\left(x_1^1,g(x_2^1,x_1^1),x_1^2\right)
=_{\{x_1^1,x^1_2,x_3^1,x_1^2,x_2^5\}}
g(x_2^1,x_3^1)
\end{equation}
\eexam

For later use, we need the following definition:

\defn\label{concdef} The \textbf{most concrete term} associated with an
expression $e$ is defined to be the unique term $t$ with the properties that
$\Expr[t]=e$ and $\Var[t]=\VarSet[e]$. The \textbf{most concrete equation}
associated with two expressions $e$ and $e'$ is defined to be the unique
equation $E$ such that $\Left[E]=e$, $\Right[E]=e'$,  and
$\Var[E]=\VarSet[e]\union\VarSet[e']$. \edefn

\exam We continue Example~\ref{eqex2}. The most concrete equation
associated with the expressions
$f\left(x_1^1,g(x_2^1,x_1^1),x_1^2\right)$
and
$g(x_2^1,x_3^1)$
is
$$f\left(x_1^1,g(x_2^1,x_1^1),x_1^2\right) =_{\{x_1^1,x_2^1,x_3^1,x_1^2\}}
g(x_2^1,x_3^1)$$
The most concrete term associated with
$f\left(x_1^1,g(x_2^1,x_1^1),x_1^2\right)$
is
$$\left(f\left(x_1^1,g(x_2^1,x_1^1),x_1^2\right),\{x_1^1,x_2^1,x_1^2\},
\sigma^5\right)$$
The most concrete term associated with
$g(x_2^1,x_3^1)$
is
$$\left(g(x_2^1,x_3^1),\{x_2^1,x_3^1\},\sigma^5\right)$$
in which we must conclude that $\Outp[g]=\sigma^5$ because it is
equated with an expression whose head is $f$.

\eexam

\section{Equational deduction}

Goguen and Meseguer~\shortcite{gogmes} prove that the following rules
for equational deduction in multisorted equational deduction
are sound and complete.

\begin{description}
\item[reflexivity] $\Frac{}{e=_Ve}$.
\item[symmetry] $\Frac{e=_Ve'}{e'=_Ve}$.
\item[transitivity] $\Frac{e=_Ve'\quad e'=_Ve''}{e=_Ve''}$.

\item[concretion]
Given a set $V$ of typed variables, $x\in V$ and an equation
$e=_Ve'$ such that $x\in
V\setminus(\VarSet[e]\union\VarSet[e'])$, and given that
$\Type[x]$ is inhabited,
$$\Frac{e=_V e'}{e=_{V\backslash \{x\}}e'}$$

\item[abstraction] Given a set $V$
of typed variables and
$x\in\Vbl[\Ssc]\setminus V$, $$\Frac{e=_V
e'}{e=_{V\union\{x\}}e'}$$

\item[substitutivity]
Given a set
$V$ of typed variables, $x\in V$, and
expressions $u$ and $u'$ for which
$\Type[x]=\Type[u]=\Type[u']$,
$$\Frac{e=_Ve'\quad u=_Wu'}{e[x\leftarrow
u]=_{V\setminus\{x\}\union W}e'[x\leftarrow u']}$$

\end{description}

\section{The sketch associated to a signature}\label{sksigg}
We now show how to construct a finite-product sketch $\S$
corresponding to a given signature in such a way that the
categories of models of the signature and of the sketch are
naturally equivalent.

Given
a signature
$\Ssc=(\Sigma,\Omega)$, we now construct a $\FinProd$ sketch
$\Sk[\Ssc]$. This sketch, like any finite-product sketch,
determines and is determined (up to isomorphism) by a
finite-product form $\F$: Precisely (see
\cite{logstr}, Section~6), there is a diagram
$\delta:I\to\FinProd$ and a global element $\Fa:1\to \v$ in
$\SynCat[\FinProd,\F]$, where $\v$ is the limit of $\delta$,
with the property that the value of $\Fa$ in the initial model
of $\SynCat[\FinProd,\F]$ in $\Set$ consists (up to isomorphism)
of the graph, diagrams and (discrete) cones that make up the
sketch $\S$.  Moreover, the finite-product theory
$\FPTh\bigl[\Sk[\Ssc]\bigr]$ (defined in \cite{ctcs},
section~7.5) is equivalent as a category to the finite-product
category $\CatTh\bigl[\FinProd,\F]\bigr]$ as defined in
\cite{logstr}.

\subsection{The graphs and cones of $\Sk[\Ssc]$}\label{sksc}
In what follows, we recursively define arrows and commutative
diagrams in $\Sk[\Ssc]$ associated to terms and equations of
$\Ssc$ respectively.

\defn The set of nodes
of $\Sk[\Ssc]$
consist by definition of the following:
\blist{OS}
\item Each type of $\Ssc$
is a node.
\item\label{listnode} Each list $v=\lstsb{\typevar}{1}{n}$
that is the input
type list
(see
Remark~\ref{inptlrem}) of at least one
operation in $\Omega$ is a node. \elist \edefn

\defn The arrows of
$\Sk[\Ssc]$ consist by definition of the following:\blist{AS}
\item Each operation $f$ in $\Omega$ is an arrow $f:
\Inp[f]\to\Outp[f]$. \item For each list
$v=\lstsb{\typevar}{1}{n}$ that is the input type list of some
operation in $\Omega$, there is an arrow
$\Proj[i]:v\to\typevar_i$ for each $i\in 1\twodots n$. (We will
write $\Proj[v,i]$ for $\Proj[i]$ if necessary to
to avoid confusion, and on the other hand we will write
$p_i$ for $\Proj[i]$ in some diagrams to save space.)\elist
\edefn

\defn\label{conedef}
The cones of $\Sk[\Ssc]$ consist by definition of the following:
For each list $v=(\typevar_1,\ldots,\typevar_n)$ that is the
input type list of some operation in $\Omega$, there is a cone
of $\Sk[\Ssc]$
with vertex $v$ and an arrow $\Proj[i]:v\to\typevar_i$ for
each $i\in 1\twodots n$.

It follows that
in a model $M$ of the sketch $\Sk[\Ssc]$,
$M(v)=\prd_{i\in 1\twodots n}M(\typevar_i)$.
\edefn

\subsubsection{Constants}
If the signature contains constants, then one of the
lists mentioned in OS.\ref{listnode} is the empty list.
As a consequence, the sketch will contain an empty cone by
Definition~\ref{conedef}, and the vertex will become a terminator in a
model.

\subsection{Terms as arrows}
\label{termsandarr} We now describe how to associate
each term of a signature $\Ssc$ to an arrow in
$\CatTh\bigl[\FinProd,\Sk[\Ssc]\,\bigr]$
and each
equation to a commutative diagram or a pair of equal arrows
in $\CatTh\bigl[\FinProd,\Sk[\Ssc]\,\bigr]$.  The
constructions given here are an elaboration of those
in \cite{ctcs}, pages~185--186.

\subsubsection{The arrow in
$\CatTh\bigl[\FinProd,\Sk[\Ssc]\,\bigr]$ corresponding to a
term}\label{conarr}
We first define recursively two arrows $Q[e]$ and $I[e]$ of
$\CatTh\bigl[\FinProd,\Sk[\Ssc]\,\bigr]$
for each expression $e$, and an arrow $D[t]$ of
$\CatTh\bigl[\FinProd,\Sk[\Ssc]\,\bigr]$ for each
term $t$.
The arrow $\Arr[t]\ceq Q\bigl[\Expr[t]\,\bigr]\o
I\bigl[\Expr[t]\,\bigr] \o D[t]$ will then be the
arrow of $\CatTh\bigl[\FinProd,\Sk[\Ssc]\,\bigr]$
associated with $t$;
the meaning of the term $t$ in a model of the signature is up to
equivalence the same function as the value of $\Arr[t]$ in the corresponding
model of $\Sk[\Ssc]$.

In
these definitions,
we suppress mention of the universal model of $\Sk[\Ssc]$.  For example,
if the universal model is
$\UnivMod[\Ssc]:\Sk[\Ssc]\to\CatTh\bigl[\FinProd,\Sk[\Ssc]\,\bigr]$ and
$\Theta$ is a node of $\Sk[\Ssc]$, then we write $\Theta$ instead of
$\UnivMod[\Theta]$.  We treat arrows of $\Sk[\Ssc]$ similarly.

\defn\label{Qcon}
For an expression $e$, $Q[e]$ is defined recursively by these
requirements:\blist{Q}
\item If $e$
is a variable $x$ of type $\tau$,
then $Q[e]\ceq\Id[\tau]$.
Using the notation introduced in A.\ref{doubleind},
if $e=x_j^i$, then
$Q[e]=\Id[\sigma^i]$,\item Suppose $e=f\lstsb{e}{1}{n}$, where $f$
is an operation with $\Inp[f]=\lstsb{\typevar}{1}{n}$ and
$\Outp[f]=\tau$.
By definition,
for $i\in(1\twodots n)$, $\Type[e_i]=\typevar_i$.

Then $Q[e]$ is defined to be the arrow
{\wider\wider\begin{equation}\label{diagtop40}
\xymatrix{ \prd_{i=1}^n\Dom[Q(e_i)] \ar[rr]^{\prd_{i=1}^n
Q[e_i]} && \crs{\typevar}{n}\ar[r]^{\,\,\,f} & \tau}
\end{equation}} \elist \edefn

\rem
We note that if $n=0$ in Q.2, in other words
$\Inp[f]$ is empty, the composite in~\eqref{diagtop40}
becomes $$\xymatrix{ 1\ar[r] & 1\ar[r]^f & \tau} $$ \erem

\defn
Let $e$ be the expression described in
Definition~\ref{Qcon} Q.2, so that
$\TypeList[e]$ is the concatenate
$\TypeList[e_1]\cdots \TypeList[e_n]$.
Then $I[e]$ is defined to be the
canonical isomorphism
$$I[e]:\prd\TypeList[e]\to\prd_{i=1}^n\TypeList[e_i]$$
given by the associative law for categorial products
in $\CatTh\bigl[\FinProd,\Sk[\Ssc]\,\bigr]$.
\edefn


\defn\label{dtdef} Let $t$ be an arbitrary term.
Then $$D[t]:\prd\InputTypes[t]\to
\prd\TypeList\bigl[\Expr[t]\,\bigr]$$ is defined to be
the unique arrow induced by requiring
that the following diagrams commute for each pair $$(i,k)\in
\bigl(1\twodots\ab{\InputTypes[t]}\bigr)\x
\bigl(1\twodots\ab{\TypeList\left[\Expr[t]\,\right]}\bigr)$$
with the property that the $i$th variable from the left in
$\Expr[t]$ is $(\Var[t])_k$.

{ \spreaddiagramcolumns{-1.5em}
\begin{equation}\label{diagbottom40}
\xymatrix{ \prd\InputTypes[t]
\ar[rr]^{D[t]}
\ar[dr]_{\Proj[k]} & &
\prd\TypeList\bigl[\Expr[t]\,\bigr]
\ar[dl]^{\Proj[i]}\\
&(\InputTypes[t])_{k} &\\
}
\end{equation}
}
\edefn

Alternatively, suppose $\Var[t]$ has length $L$ and
$\VarSet\bigl[\Expr[t]\,\bigr]$ has length $M$.
Let $\phi:1\twodots L\to 1\twodots M$ be defined by
$\phi(l)=m$ if $(\Var[t])_l=(\VarList\bigl[\Expr[t]\,\bigr])_m$
(there is a unique $m$ that makes this true).
Then we may also define
$$D[t]:\prd\InputTypes[t]\to
\prd\TypeList\bigl[\Expr[t]\,\bigr]$$ to be
the arrow $(\Proj[\phi(l)]\mid l\in 1\twodots L)$.

This works because the $\left(\phi(l)\right)$th type in
$\prd\TypeList\bigl[\Expr[t]\,\bigr]$ is indeed the type of the
$\left(\phi(l)\right)$th variable in
$\VarSet\bigl[\Expr[t]\,\bigr]_m$ (see Definition~\ref{def234}).

This is equivalent to requiring the
diagrams\eqref{diagbottom40} to commute.  The two definitions
are useful for different sorts of calculations and are therefore
included.

\exam Consider $e\ceq g(x^1_1,c)$, where $c$ is a constant of type
$\sigma^2$ and $g$ has type $\sigma^3$.  Suppose
$$t=(g(x^1_1,c),\{x^1_1,x^4_1\},\sigma^3)$$
Then $e$ corresponds to the arrow
{\taller\taller
$$
\xymatrix{ & \sigma^1\x\sigma^4
\ar[d]_{\Proj[1]}^{\brspace\bigg\} D[t]}\\ & \sigma^1
\ar[d]_{<\Id[\sigma^1],!>}^{\brspace\bigg\} I[t]}\\ & \sigma^1\x
1 \ar[d]_{\Id[\sigma^1]\xord c} \ar@{}[dd]^{\brspace\biggg\}
Q[t]} \\ & \sigma^1\x\sigma^2 \ar[d]_g\\ & \sigma^3\\ } $$ Note
that one does not have to consider constants in constructing
$D[t]$.}
\eexam

\exam
Let $e\ceq f\left(x_1^1,g(x_2^1,x_1^1),x_1^2\right)$
with $\Inp[g]=(\sigma^1,\sigma^1)$,
$\Outp[g]=\sigma^2$.
$\Inp[f]=(\sigma^1,\sigma^2,\sigma^2)$,
and $\Outp[f]=\sigma^5$.
Let
$$t\ceq(e,\{x_1^1,x_2^1,x_1^2,x_3^4\},\sigma^5)$$
Then $$\VarList\bigl[\Expr[t]\,\bigr]=(x_1^1,x_2^1,x_1^1,x_1^2)$$
$$\InputTypes[t]=(\sigma^1,\sigma^1,\sigma^2,\sigma^4)$$
$$\Var[t]=\{x_1^1,x_2^1,x_1^2,x_3^4\}$$
and
$$\TypeList\bigl[\Expr[t]\,\bigr]=(\sigma^1,\sigma^1,\sigma^1,\sigma^2)$$

If we use the first definition of $D[t]$ in Definition~\ref{dtdef},
then the following four triangles must commute:

{ \spreaddiagramcolumns{-1.5em}
\mathchardef\times="0202

$$ \begin{array}{@{\hspace{-1.5em}}cc}
\xymatrix{ \sigma^1\x\sigma^1\x\sigma^2\x\sigma^4 \ar[rr]^{D[t]}
\ar[dr]_{\Proj[1]} & &                                      x
\sigma^1\x\sigma^1\x\sigma^1\x\sigma^2
\ar[dl]^{\Proj[1]}\\
&\sigma^1&\\
}
 &
\xymatrix{ \sigma^1\x\sigma^1\x\sigma^2\x\sigma^4 \ar[rr]^{D[t]}
\ar[dr]_{\Proj[2]} & &
\sigma^1\x\sigma^1\x\sigma^1\x\sigma^2
\ar[dl]^{\Proj[2]}\\
&\sigma^1&\\
}
\end{array}$$
\begin{equation}\label{page41}
 \begin{array}{@{\hspace{-1.5em}}cc}
\xymatrix{ \sigma^1\x\sigma^1\x\sigma^2\x\sigma^4 \ar[rr]^{D[t]}
\ar[dr]_{\Proj[1]} & &
\sigma^1\x\sigma^1\x\sigma^1\x\sigma^2
\ar[dl]^{\Proj[3]}\\
&\sigma^1&\\
}
&
\xymatrix{ \sigma^1\x\sigma^1\x\sigma^2\x\sigma^4 \ar[rr]^{D[t]}
\ar[dr]_{\Proj[3]} & &
\sigma^1\x\sigma^1\x\sigma^1\x\sigma^2
\ar[dl]^{\Proj[4]}\\
&\sigma^2&\\
}
 \end{array}
\end{equation}
}
It follows that $D[t]$ is given by the following diagram, where
to save space we write $p_k$ for $\Proj[k]$.
\begin{equation}\label{pppp41} \xymatrix{
\sigma^1\xord\sigma^1\xord\sigma^2\xord\sigma^4
\ar[rrr]^{<p_1,p_2,p_1,p_3>} &&&
\sigma^1\xord\sigma^1\xord\sigma^1\xord\sigma^2\\
}
\end{equation}
and that $\Arr[t]$ is   the composite
$$\xymatrix{
\sigma^1\xord\sigma^1\xord\sigma^2\xord\sigma^4
\ar[d]_{<p_1,p_2,p_1,p_3>}^{\brspace\bigg\}
D[t]} \\ \sigma^1\xord\sigma^1\xord\sigma^1\xord\sigma^2
\ar[d]_{<p_1,<p_2,p_3>,p_4>}^{\brspace\bigg\}
I[t]} \\ \sigma^1\xord(\sigma^1\xord\sigma^1)\xord\sigma^2
\ar[d]_{\Id[\sigma^1]\xord g\xord\Id[\sigma^2]}
\ar@{}[dd]^{\brspace\biggg\} Q[t]}
\\
\sigma^1\xord\sigma^2\xord\sigma^2
\ar[d]_f
\\
\sigma^5
}$$\eexam

\subsection{The diagram associated to an equation}

Let the equation $E\ceq e=_Ve'$ be given.
Define the terms $t_1$ and $t_2$ by
$t_1=(e,\Var[E],\Type[E])$ and $t_2=(e',\Var[E],\Type[E])$
(using the
notation of~\ref{termdeff}).
Recall that $\Type[E]=\Type[e]=\Type[e']$.
The notation $\InputTypes[E]$ will denote the list
$\InputTypes[t_1]$, which is the same as $\InputTypes[t_2]$.
As in~\ref{conarr}, we have
arrows $\Arr[t_1]$ and $\Arr[t_2]$ with the same domain and codomain.
We will associate the diagram
\begin{equation}\label{p42first}
\xymatrix{
\InputTypes[E]\ar[rr]<1ex>^{\quad\Arr[t_1]}
\ar[rr]<-1ex>_{\quad\Arr[t_2]} &&
\Type[E]\\
}
\end{equation}
to the equation $E$.  By~\ref{conarr}, this is the same as
{\taller
\begin{equation}\label{p42second}
\xymatrix{
\InputTypes[E]\ar[rr]^{\quad D[t_1]} \ar[d]|{D[t_2]} &&
\TypeList[e]\ar[d]|{Q[t_1]\o
I[t_1]}\\
\TypeList[e'] \ar[rr]_{\quad Q[t_2]\o I[t_2]} && \Type[E] \\
}
\end{equation}
}

This completes the translation.

\rem The commutative diagram as exhibited above can also be
viewed as a pair of formally
equal arrows as in Diagram~\eqref{p42first},
and in what follows we will use this description frequently.\erem

\subsection{Examples}\label{twoex}
We work out below two examples in detail to facilitate later discussion of
substitution.

\exam\label{fexam} $$ \begin{array}{l}
e \ceq f\left( \ds x11, \ds x43, \ds x32, \ds x11, g( \ds x11, \ds x32), \ds
x21\right)\\ \\
\begin{cases}\Inp[f] = \sigma^1\x \sigma^4\x \sigma^3\x \sigma^1\x \sigma^5\x
\sigma^2 & \\ \Outp[f] = \sigma^5\end{cases}\\ \\
\begin{cases}\Inp[g] =
\sigma^1\x \sigma^3&\\ \Outp[g]= \sigma^5&\end{cases}\\
\\ V=\{\ds x11,\ds
x12,\dsu x13,\ds x21,\dsu x22,\dsu x31,\ds x32,\ds x43\} \end{array}$$
The
underlined variables are redundant; that is, they do not appear in the
expression $e$.
$$ \begin{array}{l}
\begin{cases}u \ceq h( \ds x21, \ds
x32)&\\ \Inp[h] = \sigma^2\x \sigma^3&\\ \Outp[h] = \sigma^3&\end{cases}\\ \\
W \ceq\{ \dsu x11, \ds x21, \dsu x22, \dsu x23, \dsu x31, \ds
x32, \dsu x33\} \end{array}$$ $ \ds x32$ is a variable for which
we are making a substitution. We wish to calculate
\begin{align*}(e,V, \sigma^5)\left[ \ds x32 \leftarrow
(u,W, \sigma^3)\right] &= \bigl(e[ \ds x32\leftarrow
u],(V\setminus\{ \ds x32\})\union W, \sigma^5\bigr)\end{align*}
By direct calculation,
$$(V\setminus\{ \ds x32\})\union W=
\{ \ds x11, \ds x12, \dsu x13, \ds x21, \dsu  x22, \dsu x23, \dsu x31,
\ds x32, \dsu x33, \ds x43\}$$
and
$$e( \ds x32\leftarrow u)=
f\left( \ds x11, \ds x43, e( \ds x21, \ds  x32), \ds x11,
g( \ds x12, e( \ds x21, \ds x32)), \ds x21\right)$$
We now exhibit the arrows for $e$ and $u$ over $V$:
$e = f\left( \ds x11, \ds x43, \ds x32, \ds x11, g( \ds x11, \ds x32), \ds
x21\right)$:
{\def\brspace{\hspace{5em}}
$$\xymatrix{
\sigma^1\xord\sigma^1\xord\mk{\sigma}^1\xord\sigma^2
\xord\mk{\sigma}^2\xord\mk{\sigma}^3\xord\sigma^3
\xord\sigma^4
\ar[d]_{<p_1,p_8,p_7,p_1,p_2,p_7,p_4>}^{\brspace\bigg\} D[e]}
\\
\sigma^1\xord\sigma^4\xord\sigma^3\xord\sigma^1 \xord\sigma^1\xord\sigma^3
\xord\sigma^2
\ar[d]_{<p_1,p_2,p_3,p_4,<p_5,p_6>,p_7>}^{\brspace\bigg\} I[e]}
\\
\sigma^1\xord\sigma^4\xord\sigma^3\xord\sigma^1\xord(\sigma^1\xord\sigma^3)
\xord\sigma^2
\ar[d]_{\Id[\sigma^1]\xord \Id[\sigma^4]\xord\Id[\sigma^3]\xord
\Id[\sigma^1]\xord g\xord\Id[\sigma^2]}
\ar@{}[dd]^{\brspace\biggg\} Q[e]}
\\
\sigma^1\xord\sigma^4\xord\sigma^3\xord\sigma^1\xord\sigma^5
\xord\sigma^2
\ar[d]_f
\\
\sigma^5
}$$
$u=h(\ds x21, \ds x32)$ (over $W$):}
$$\xymatrix{
\mk{\sigma}^1\xord\sigma^2\xord\mk{\sigma}^2
\xord\mk{\sigma}^2\xord\mk{\sigma}^3\xord\sigma^3
\xord\mk{\sigma}^3
\ar[d]_{<p_2,p_6>}^{\brspace\bigg\} D[u]}\\
\sigma^2\xord\sigma^3
\ar[d]_{<p_1,p_2>=\Id[\sigma^2\xord\sigma^3]
=\Id[\sigma^2]\xord\Id[\sigma^3]}^{\brspace\bigg\} I[u]}\\
\sigma^2\xord\sigma^3
\ar[d]_h^{\brspace\bigg\} Q[u]}\\
\sigma^3}$$

Therefore
$\Arr[u,W,\sigma^3]=h<p_1,p_2><p_2,p_6>=h<p_2,p_6>$
and $u=h( \ds x21, \ds x32)$ (over $\left(V\setminus\{ \ds
x32\}\right)\union W$) is the arrow
$$\xymatrix{
\mk{\sigma}^1\xord\mk{\sigma}^1\xord\mk{\sigma}^1
\xord\sigma^2\xord\mk{\sigma}^2\xord\mk{\sigma}^2
\xord\mk{\sigma}^3\xord\sigma^3\xord\mk{\sigma}^3
\xord\mk{\sigma}^4
\ar[d]_{<p_4,p_8>}^{\brspace\bigg\} D[u]}\\
\sigma^2\xord\sigma^3
\ar[d]_{<p_1,p_2>=\Id[\sigma^2\xord\sigma^3]
=\Id[\sigma^2]\xord\Id[\sigma^3]}^{\brspace\bigg\}I[u]}\\
\sigma^2\xord\sigma^3
\ar[d]_h^{\brspace\bigg\}Q[u]}\\
\sigma^3}$$
so that
$\Arr[u,\left(V\setminus\{ \ds x32\}\right)\union W,\sigma^3]
=h<p_1,p_2><p_4,p_8>=h<p_4,p_8>$.

Note\label{twoexx}
that we have the maps
\begin{gather*}
\alpha:\Type[W]\to\Type\bigl[\left(V\setminus\{ \ds
x32\}\right)\union W\bigr]\\ \beta:\Type\bigl[\left(V\setminus\{
\ds x32\}\right)\union W\bigr]\to\Type[W] \end{gather*} as shown
below:
{\taller\taller$$\xymatrix{
\Type[W]=\mk{\sigma}^1\xord\sigma^2\xord\mk{\sigma}^2
\xord\mk{\sigma}^2\xord\mk{\sigma}^3\xord\sigma^3
\xord\mk{\sigma}^3
 \ar[d]|{\alpha=< \ds x11!, \ds x12!, \ds x13!, p_2, \ds x22!, \ds x32!, \ds
x31!, p_6, \ds x33!, \ds x43!>}\\
\mk{\sigma}^1\xord\mk{\sigma}^1
\xord\mk{\sigma}^1\xord\sigma^2\xord\mk{\sigma}^2\xord\mk{\sigma}^2
\xord\mk{\sigma}^3\xord\sigma^3\xord\mk{\sigma}^3\xord\mk{\sigma}^4
}$$
where
the codomain is
$$\Type\bigl[\left(V\setminus\{ \ds x32\}\right)\union W\bigr]$$
and where we have identified the variable $x^i_j$ with $$
\xymatrix{ W \ar[r]^{!} & 1 \ar[r]^{x^i_j} & \sigma^i }$$ The
map $\beta $ is similarly defined.
It follows that
$$<p_4,p_8>\alpha=<p_2,p_6>$$ and $$<p_2,p_6>\beta=<p_4,p_8>$$
and that $$\Arr[u,W,\sigma^3]=\Arr\bigl[u,\left(V\setminus\{ \ds
x32\}\right)\union W, \sigma^3\bigr]\o\alpha$$
and
$$\Arr\bigl[u,\left(V\setminus\{ \ds x32\}\right)\union W,
\sigma^3\bigr]=\Arr[u,W,\sigma^3]\o\beta$$
These observations,
although made in this special case, capture general features of
the system that we shall need later.  We record these in passing
in the following lemmas.  The proofs are quite straightforwarded
and are omitted in view of the preceding example.
\thl\label{lemmaa} Let
$t=\left(\Expr[t],\Var[t],\Type[t]\right)$ be any term. Then
$$Q[t]=Q\bigl[\Expr[t]\bigr]=Q\bigl[\Expr[t],\Var[\Expr[t]],\Type[t]\bigr]$$
and
$$I[t]=
I\bigl[\Expr[t]\bigr]=I\bigl[\Expr[t],\Var[\Expr[t]],\Type[t]\bigr]$$
are determined by the most concrete term associated with $t$ as defined
in~\ref{concdef}.  These specifically do not depend on
$\Var[t]$. \ethl

Note that by contrast $D[t]$ does depend on $\Var[t]$.

\thl\label{lemmb} Let $e$ be an expression of type $\tau$, and let $V_1$ and
$V_2$ be lists of variables such that $\VarSet[e]\includedin V_1$ and
$\VarSet[e]\includedin V_2$. Let
$$t_1\ceq\left[e,V_1,\tau\right]$$ and
$$t_2\ceq\left[e,V_2,\tau\right]$$
Then there are arrows
$$\alpha_{12}:\prd\TypeList[t_1]\to\prd\TypeList[t_2]$$
and
$$\alpha_{21}:\prd\TypeList[t_2]\to\prd\TypeList[t_1]$$
for which
$\Arr[t_1]=\Arr[t_2]\o\alpha_{12}$ and
$\Arr[t_2]=\Arr[t_1]\o\alpha_{21}$.
\ethl

\rem
Lemmas~\ref{lemmaa} and~\ref{lemmb} will be used later in our discussion of
the rules
``concretion'' and ``abstraction'' that have to do with including extraneous
variables in and excluding them from the list of variables of some term.
\erem

We now proceed with our example.  After
substitution, $$e(\ds x32 \leftarrow
u)=f\left(\ds x11, \ds x43, h( \ds x21, \ds x32), \ds x11, g(
\ds x12, h( \ds x21, \ds x32)), \ds x21\right)$$ This
corresponds to the arrow {\def\brspace{\hspace{7em}} $$
\xymatrix{ \sigma^1 \xord \sigma^1 \xord \mk{\sigma}^1 \xord
\sigma^2 \xord \mk{\sigma}^2 \xord \mk{\sigma}^2 \xord
\mk{\sigma}^3 \xord \sigma^3 \xord \mk{\sigma}^3 \xord \sigma^4
\ar[d]_{<p_1, p_{10}, p_4, p_8, p_1, p_2, p_4, p_8, p_4>}^{\brspace\bigg\}D[e]}\\
\sigma^1 \xord \sigma^4 \xord \sigma^2 \xord
\sigma^3 \xord \sigma^1 \xord \sigma^2 \xord \sigma^3 \xord \sigma^2
\ar[d]_{\left<p_1, p_2, <p_3, p_4>, \left<p_4, <p_7, p_8>\right>,
p_9\right>}^{\brspace\bigg\}I[e]}\\
\sigma^1 \xord \sigma^4 \xord (\sigma^2 \xord
\sigma^3 ) \xord \sigma^1 \xord \left(\sigma^1 \xord (\sigma^2 \xord
\sigma^3)\right) \xord \sigma^2
\ar[d]_{\Id[\sigma^1] \xord \Id[\sigma^4] \xord
\Id[\sigma^2 \xord \sigma^3] \xord \Id[\sigma^1] \xord
\left(\Id[\sigma^1] \xord h\right) \xord \Id[\sigma^2]}
\ar@{}[ddd]^{\brspace\hspace{1em}\Biggggg\} Q[e]}\\
\sigma^1 \xord \sigma^4 \xord (\sigma^2 \xord
\sigma^3 ) \xord \sigma^1 \xord (\sigma^1 \xord \sigma^3)
\xord \sigma^2
\ar[d]_{\Id[\sigma^1] \xord \Id[\sigma^4] \xord h \xord
\Id[\sigma^1] \xord g \xord \Id[\sigma^2]}\\
\sigma^1 \xord \sigma^4 \xord \mk{\sigma}^3 \xord
\sigma^1 \xord \mk{\sigma}^5 \xord \sigma^2
\ar[d]_f\\
\sigma^5
}$$
}

We now calculate
\begin{equation*}\begin{split}
e\left[x^3_2\leftarrow u\right]
&=
f\left(\Id[\sigma^1] \xord \Id[\sigma^4] \xord e \xord
\Id[\sigma^1] \xord g \xord \Id[\sigma^3]\right)\\
& \quad\quad \o \left( \Id[\sigma^1] \xord \Id[\sigma^4] \xord
\Id[\sigma^2 \xord \sigma^3] \xord \Id[\sigma^1] \xord
( \Id[\sigma^1] \xord e)
\xord \Id[\sigma^2]\right) \\
& \quad\quad \o < p_1, p_2, < p_3, p_4 >, p_5, < p_6, < p_7, p_8>>, p_9>\\
& \quad\quad \o
< p_1, p_{10}, p_4, p_8, p_1, p_2, p_4, p_8, p_4>\\
&=
f\bigl( \Id[\sigma^1] \xord \Id[\sigma^4] \xord
\Id[\sigma^3] \xord \Id[\sigma^1] \xord g \xord
\Id[\sigma^2]
\bigr)\\
& \quad\quad \o
< p_1, p_2, e<p_3, p_4>, p_8, <p_6, e<p_7, p_8>>, p_9>\\
& \quad\quad \o
<p_1, p_{10}, p_4, p_8, p_1, p_2, p_4, p_8, p_4>\\
&= f\bigl(
\Id[\sigma^1] \xord \Id[\sigma^4] \xord
\Id[\sigma^3] \xord \Id[\sigma^1] \xord g \xord \Id[\sigma^2]
\bigr)\\
& \quad\quad \o
< p_1 , p_2 , p_3 , p_4 , < p_5 , p_6> , p_7>\\
& \quad\quad \o
\left<p_1 , p_2, e<p_3, p_4>, p_5, p_6, e<p_7, p_8>, p_9\right>\\
& \quad\quad \o
<p_1, p_{10}, p_4, p_8, p_1, p_2, p_4, p_8, p_4>\\
&=
f\bigl(
\Id[\sigma^1] \xord \Id[\sigma^4] \xord
\Id[\sigma^3] \xord \Id[\sigma^1] \xord g \xord \Id[\sigma^2]
\bigr)\\
& \quad\quad \o
\left<p_1, p_2, p_3, p_4, <p_5, p_6>, p_7\right>\\
& \quad\quad \o
\left<p_1, p_{10}, e<p_4, p_8>, p_1, p_2, e<p_4, p_8>, p_4\right>\\
&=
f\bigl(
\Id[\sigma^1] \xord \Id[\sigma^4] \xord
\Id[\sigma^3] \xord \Id[\sigma^1] \xord g \xord \Id[\sigma^2]
\bigr)\\
& \quad\quad \o \left< p_1, p_2, p_3, p_4, <p_5, p_6>, p_7\right>\\
& \quad\quad \o
<p_1, p_{10}, p_8, p_1, p_2, p_8, p_4>\\
& \quad\quad \o <p_1, p_2, p_3, p_4, p_5, p_6, p_7, e<p_4,p_8>, p_9, p_{10}>\\
&=\Arr\left[e,V\union W, \Type[e]\right]\\
& \quad\quad \o \left<p_1, p_2, p_3, p_4, p_5, p_6, p_7,
\Arr\left[u,\left(V\setminus\{x^3_2\}\right)\union W,
\Type[u]\right],p_9,p_{10}\right>
\end{split}\end{equation*}

\exam\label{sexam}

$$u=m( x^2_1, x^2_1, x^4_4)$$
$$\begin{cases} \Inp[u]= {\sigma}^2 \xord \sigma^2 \xord \sigma^4 &\\
\Outp[u] = \sigma^3 &\\
\end{cases}$$
$$W=\{ \mk{x}^1_1, x^2_1, \mk{x}^2_2, \mk{x}^2_3, \mk{x}^3_1, \mk{x}^3_3,
x^4_4\}$$ This is different from Example~\ref{fexam} because the variable
$x^3_2$ (in $e$) for which we are making the substitution does not reappear in
$u$.

This is the arrow for $u=m(x^2_1, x^2_1, x^4_4)$ over $W$:
$$\xymatrix{
\sigma^1 \xord \sigma^2 \xord
\mk{\sigma}^2 \xord \mk{\sigma}^2 \xord
\mk{\sigma}^3 \xord \mk{\sigma}^3 \xord \sigma^4
\ar[d]_{<p_1, p_1, p_7>}^{\brspace\bigg\} D[u]}\\
\sigma^2 \xord \sigma^2 \xord \sigma^4
\ar[d]_{<p_1, p_2, p_3>
=\Id[\sigma^2 \xord \sigma^2 \xord \sigma^4]
= \Id[\sigma^1] \xord \Id[\sigma^2] \xord \Id[\sigma^4]}^{\brspace\bigg\} I[u]}\\
\sigma^2 \xord \sigma^2 \xord \sigma^4
\ar[d]_m^{\brspace\bigg\} Q[u]}\\
\sigma^5
}$$
We calculate
$$\left(V\setminus\{x^3_2\}\right)\union W=
\{x^1_1, x^1_2, \mk{x}^1_3, x^2_1, \mk{x}^2_2, \mk{x}^2_3, \mk{x}^3_1,
\mk{x}^3_3, x^4_3, x^4_4\}$$
Then $u=m(x^2_1, x^2_1, x^4_4)$ (over
$\bigl(V\setminus\{x^3_2\}\bigr)\union W$) gives the arrow
$$\xymatrix{ \sigma^1 \xord \sigma^1 \xord \mk{\sigma}^1 \xord
\sigma^2 \xord \mk{\sigma}^2 \xord \mk{\sigma}^2 \xord
\mk{\sigma}^3 \xord \mk{\sigma}^3 \xord \sigma^4  \xord \sigma^4
\ar[d]_{<p_4, p_4, p_{10}>}^{\brspace\bigg\} D[u]}\\ \sigma^2
\xord \sigma^2 \xord \sigma^4 \ar[d]_{<p_1, p_2, p_3>
=\Id[\sigma^2 \xord \sigma^2 \xord \sigma^4] = \Id[\sigma^1]
\xord \Id[\sigma^2] \xord \Id[\sigma^4]}^{\brspace\bigg\}
I[u]}\\ \sigma^2 \xord \sigma^2 \xord \sigma^4
\ar[d]_m^{\brspace\bigg\} Q[u]}\\ \sigma^5 }$$

After substitution,
$$e[x^3_2\leftarrow u]=
f\bigl( x^1_1, x^4_3, m\left( x^2_1, x^2_1, x^4_4\right), x^1_1,
g\left(x^1_2, m\left(x^2_1, x^2_1, x^4_4\right)\right),
x^2_1\bigr)$$ This corresponds to the arrow shown below.
{\taller\taller\def\brspace{\hspace{10em}}$$ \xymatrix{
\sigma^1 \xord \sigma^1 \xord \mk{\sigma}^1 \xord \sigma^2 \xord
\mk{\sigma}^2 \xord \mk{\sigma}^2 \xord \mk{\sigma}^3 \xord \mk{\sigma}^3
\xord \sigma^4 \xord \sigma^4 \ar[d]|{\left<p_1, p_9, p_{10}, p_4,
p_4, p_{10}, p_1, p_2, p_4, p_4, p_{10},
p_2\right>}^{\brspace\bigg\} D[e]}\\ \sigma^1 \xord \sigma^4 \xord
\sigma^2 \xord \sigma^2 \xord \sigma^4 \xord \sigma^1 \xord \sigma^1 \xord
\sigma^2 \xord \sigma^2 \xord \sigma^4 \xord \sigma^2 \ar[d]|{\bigl<p_1,
p_2, \left<p_3, p_4, p_5\right>, p_6, \left<p_7, <p_8, p_9,
p_{10}>\right>, p_{11}\bigr>}^{\brspace\bigg\} I[e]}\\ \sigma^1
\xord \sigma^4 \xord \left( \sigma^2 \xord \sigma^2 \xord \sigma^2 \right)
\xord \sigma^1 \xord \bigl( \sigma^1 \left( \sigma^2 \xord \sigma^2 \xord
\sigma^4\right)\bigr) \xord \sigma^2 \ar[d]|{\Id[\sigma^1] \xord
\Id[\sigma^4] \xord \Id\left[\sigma^2 \xord \sigma^2 \xord
\sigma^4\right] \xord \Id[\sigma^2] \left( \Id[\sigma^1] \xord h
\right) \xord \Id[\sigma^2]}\\ \sigma^1 \xord \sigma^4 \xord \left(
\sigma^2 \xord \sigma^2 \xord \sigma^2 \right) \xord \sigma^1 \xord \bigl(
\sigma^1 \xord \mk{\sigma}^3 \bigr) \xord \sigma^2
\ar[d]|{\Id[\sigma^1] \xord \Id[\sigma^4] \xord h \xord
\Id[\sigma^1] \xord g \xord \Id[\sigma^2]}^{\brspace\biggggg\}
Q[e]}\\ \sigma^1 \xord \sigma^4 \xord \mk{\sigma}^3 \xord
\sigma^1 \xord \mk{\sigma}^5 \xord \sigma^2 \ar[d]^f\\ \tau }$$}

We next re-express this in a convenient form:

\begin{equation*}\begin{split}
e[x^3_2\leftarrow u]&=
f \o \bigl( \Id[\sigma^1] \xord \Id[\sigma^4] \xord h \xord
\Id[\sigma^2] \xord g \xord \Id[\sigma^2] \bigr)\\
 & \quad\quad \o \bigl(
\sigma^1 \xord \sigma^4 \xord \Id\left[ \sigma^2
\xord \sigma^2 \xord \sigma^4\right] \xord \Id[\sigma^1] \xord
\left(  \Id[\sigma^1] \xord h \right) \xord \Id[\sigma^2] \bigr)\\
& \quad\quad \o \bigl< p_1, p_2, \left< p_3, p_4, p_5 \right>,
p_6, \left< p_7, \left< p_8, p_9, p_{10} \right> \right>,
p_{11}\bigr>\\
 & \quad\quad \o \bigl< p_1, p_9, p_4, p_4, p_{10},
p_1, p_2, p_4, p_4, p_{10}, p_2 \bigr>\\
 &=f \o \bigl(
\Id[\sigma^1] \xord \Id[\sigma^4] \xord \Id[\sigma^3] \xord
\Id[\sigma^1] \xord g \xord \Id[\sigma^2] \bigr)\\
 & \quad\quad \o
\bigl< p_1, p_2, h<p_3, p_4, p_5>, p_6, \left< p_7, h<p_8, p_9,
p_{10}> \right>, p_{11} \bigr>\\
 & \quad\quad \o \bigl< p_1, p_9,
p_4, p_4, p_{10}, p_1, p_2, p_4, p_4, p_{10}, p_2 \bigr>\\
 &= f
\o \bigl( \Id[\sigma^1] \xord \Id[\sigma^4] \xord \Id[\sigma^3]
\xord \Id[\sigma^1] \xord g \xord \Id[\sigma^2] \bigr)\\
& \quad\quad \o \bigl< p_1, p_9, h<p_4, p_4, p_{10}>, p_1, \left<
p_2, h<p_4, p_4, p_{10}> \right>, p_2 \bigr> \\
 &= f \o \bigl(
\Id[\sigma^1] \xord \Id[\sigma^4] \xord \Id[\sigma^3] \xord
\Id[\sigma^1] \xord g \xord \Id[\sigma^2] \bigr)\\
 & \quad\quad \o
\bigl< p_1, p_2, p_3, p_4, \left<p_5, p_6 \right>, p_7 \bigr>\\
& \quad\quad \o \bigl< p_1, p_{10}, p_8, p_1, p_2,  p_8, p_4
\bigr>\\
 & \quad\quad\o \bigl< p_1, p_2, p_3, p_4, p_5, p_6, p_7,
h<p_4, p_4, p_{10}>, p_9, p_1, p_{11} \bigr>\\
&= \Arr \left( e,
V \union W, \sigma^5 \right)\\
 & \quad\quad \o < p_1, p_2, p_3,
p_4, p_5, p_6, p_7, \Arr \bigl(u, \left( V\setminus\{x^3_2\}
\right) \union W \bigr), p_9, p_{10}, p_{11} \bigr>
\end{split}\end{equation*}

\rem In Examples~\ref{fexam} and~\ref{sexam} we may define a map
$$A(e,u):\prd\TypeList\bigl[\left(V\setminus\{x^3_2\}\right)\union W\bigr]\to
\prd\TypeList\left[V\union W\right]$$
as follows.
Choose $I\in 1\twodots\Length \left[V\union W\right]$ such that
$$\bigl(\TypeList[V\union W]\bigr)_I=x^3_2$$
We next define
for all $ i\in1\twodots\Length\left[V\union W\right]$
$$\left(A(e,u) \right)_i=\Proj[i]$$ and $$\left(A(e,u)
\right)_I =\Arr\bigl(u, \left(V\setminus\{x^3_2\} \right) \union
W,\sigma^5\bigr)$$ With this definition in the previous two
examples we have $$\Arr\bigl[ \left(e,V,\sigma^5
\right)\left[x^3_2\leftarrow \left(u, W, \sigma^3
\right)\right]\bigr] =\Arr\left[e,V\union W,\sigma^5\right]\o
A(e,u)$$ This is again a general feature that the following
discussion of substitution is intended to capture.\erem

\subsection{Substitution}

As terms are defined recursively, substitution may be defined either by
structural recursion or, in view of Examples~\ref{fexam}
and~\ref{sexam}, using composition.  These two ways of defining
substitution are convenient for different purposes.  Here we
establish the equivalence of the two procedures.

\subsubsection{Recursive definition}
Suppose $ \left(u,W,\tau \right) $ is to be substituted for $x$ in $t=
\left(e,V,\sigma \right)$  where we assume that $x\in V$ and $\Type[x]=\tau$.
We may define this recursively as follows:
If $t= \left(x,V,\tau \right)$, then
$$t\bigl[x\leftarrow \left(u,W,\tau \right)\bigr]
\ceq
 \bigl(u, \left(V\setminus\{x\} \right)\union W,
\tau\bigr)$$
If $t=
\left(f(e_1,\ldots,e_n),V,\tau\right)$, then
\begin{equation*}
t\bigl[x\leftarrow \left(u,W,\tau \right)\bigr]
\ceq \bigl( f \left(e_1[x\leftarrow u],\ldots,e_n[x\leftarrow u]
\right), \left(V\setminus\{x\} \right)\union
W,\sigma\bigr)\end{equation*}

Note that the added complication here is
owing to the fact that the list of variables can be independently specified
and that we have tacitly assumed the usual recursive definition of the
substitution of one expression in the place of a free variable in
another. The above serves as a basis for the
recursive definition of the arrow corresponding to $t$ once the substitution
has been made.

If $t= \left(x,V,\tau \right)$, then
$$\Arr\bigl[x\leftarrow \left(u,W,\tau \right)\bigr]
\ceq
 \Arr\bigl[u, \left(V\setminus\{x\} \right)\union W,
\tau\bigr]$$
If $t=
\left(f(e_1,\ldots,e_n),V,\tau\right)$, then
\begin{equation*}
\Arr\bigl[x\leftarrow \left(u,W,\tau \right)\bigr]
\ceq \Arr\bigl[ f \left(e_1[x\leftarrow
u],\ldots,e_n[x\leftarrow u] \right), \left(V\setminus\{x\}
\right)\union W,\sigma\bigr]\end{equation*}

\subsubsection{Direct definition}
The alternative way suggested by Examples~\ref{fexam} and~\ref{sexam} is to
define $\Arr \left[e[x\leftarrow u] , \left(V\setminus\{x\} \right),
\sigma\right]$ directly, given
\begin{gather}
\Arr \left[e,V,\sigma \right]=Q[e]I[e]D[(e,V,\sigma)]\\
\Arr \left[u,W,\tau \right]=Q[u]I[u]D[(u,W,\tau)]
\end{gather}
In view of Lemma~\ref{lemmaa}, we may suppose $\Arr\left[e,V\union
W,\sigma\right]= Q[e]I[e]D[(e,V\union W,\sigma)]$, that is that
$\Arr\left[e,V\union W,\sigma\right]$ differs from
$\Arr\left[e,V,\sigma\right]$ only in the $D$-composand of the arrow.

We have
$$\Dom\bigl[D \left[(e,V\union W,\sigma)  \right]  \bigr]
= \prd \left(V\union W \right)$$
Choose $I\in 1\twodots\Length[V\union W]$ such that $(V\union W)_I=x$.
Define an arrow
$$A\ceq A\bigl( \left(e,V,\sigma \right), \left(u,W,\tau \right)\bigr):
\prd\bigl( \left(V\setminus\{x\} \right)\union W\bigr)\to
\prd \left(V\union W \right)$$
as follows:
For all $i\in \bigl(1\twodots\Length \left(V\union W
\right)\bigr)\setminus\{I\}$, $A_i=\Proj[i]$
and $$A_I=\Arr \left[u, \left(V\setminus\{x\} \right) \union W,
\tau \right]$$ Note that $\prd \bigl( \left(V\setminus\{x\}
\right)\union W\bigr)$ and $\prd \left(V\union  W\right)$ can
differ in at most one factor depending on whether $x\in W$ or
not.

Finally, we define
$$\Arr \bigl[e[x\leftarrow u],
\left(V\setminus\{x\}\right)\union W, \tau \bigr] =
\Arr \left[e, V\union W, \sigma\right]\o A\bigl( \left(
e,V,\sigma\right), \left(u,W,\tau \right)\bigr)$$
We have given
two methods of obtaining the arrow corresponding
to the term for which substitution has been made.  It remains to
be seen that these two methods give the same arrow.

\subsubsection*{Proof by structural induction}
\paragraph*{Base case:} $t= \left(x,V,\sigma \right)$ and $\sigma=\tau$.
We note that
\begin{align*}
\Arr[t] &=
Q[x]\o I[x]\o D[t]\\
&= \Id[\sigma]\o\Id[\sigma]\o \proj[I]\quad\text{(where $
\left(V\union W \right)_I=x$)}\\ &=\proj[I] \end{align*} From
the direct definition we have \begin{multline*} \Arr
\bigl[x[x\leftarrow u]  , \left(V\setminus\{x\} \right)\union
W,\sigma\bigr]
= \Arr \left[x, V\union W,\tau  \right]\o A \bigl( \left(x,V,\tau \right),
\left(u,W,\tau \right)\bigr)\\
=\proj[I]\o A \bigl( \left(x,V,\tau \right), \left(u,W,\tau
\right)\bigr) =\Arr \bigl(u, \left(V\setminus\{x\} \right)\union
W,\tau\bigr) \end{multline*} by the definition of $A \bigl(
\left(x,V,\tau \right), \left(u,W,\tau \right)\bigr)$, which
agrees with the recursive definition.
\paragraph*{Induction
step:} $t= \left(f(e_1,\ldots,e_n), V,\sigma \right)$, where
$\Outp[f]=\sigma$ and
for all $i\in 1\twodots
n$,  $\Outp[e_i]=\gamma_i$.

We note that, if we define
$$t'\ceq \left(f(e_1,\ldots,e_n) ,V',\sigma\right)$$
where $\VarSet \left[\Expr[t]  \right]\includedin V'$, then
\begin{align*}
\Arr[t'] &=
f\o \left(\prd_{1\in 1\twodots n}Q[e_i]
\right)\o\left<I[e_1],\ldots,I[e_n]\right>
\o D[t']\\
&=f\o\bigl<
Q[e_1]I[e_1]D[t'],\ldots,Q[e_n]I[e_n]D[t']
\bigr>\\
&=f\o\bigl<
\Arr \left[e_1,V',\gamma_1  \right],\ldots,
\Arr \left[e_n,V',\gamma_n  \right]\bigr>
\end{align*}
Although this last equality is obvious, a complete proof may
require a lemma.

By induction hypothesis, we have, for all $i\in 1\twodots n$,
\begin{equation*}\Arr \bigl[
e_i[x\leftarrow u], \left(V\setminus\{x\} \right)\union W,\gamma_i
\bigr]
=\Arr \bigl[
e_i,V\union W,\gamma_i)\bigr]\o A \bigl( \left(e_i,V,\gamma_i \right),
\left(u,W,\tau \right) \bigr)\end{equation*}
We next note that the direct definition yields
\begin{multline*}
\Arr[t] =\Arr\bigl[
f(e_1,\ldots,e_n)[x\leftarrow u], \left(V\setminus\{x\} \right),\sigma
\bigr]\\
= \Arr \left[f(u_1\ldots,u_n), V\union W,\sigma  \right]\o
A\bigl( \left(e,V,\sigma \right), \left(u,W,\tau \right)
\bigr)\\ \typeout{Here is A}
\shoveright{\text{where $e\ceq f(\lst e1n)$}}\\
= f\o \bigl<
\Arr \left[e_1,V\union W,\gamma_1  \right],\ldots,
\Arr \left[e_n,V\union W,\gamma_n  \right]
\bigr>\\
\quad\quad\o A\bigl(
\left(e,V,\sigma \right), \left(u,W,\tau \right)
\bigr)
\\ \typeout{Here is B}
=f\o \bigl<
\Arr \left[e_1,V\union W,\gamma_1  \right]\o
A\bigl(
\left(e,V,\sigma \right), \left(u,W,\tau \right)
\bigr),\\
\quad\quad\ldots,
\Arr \left[e_n,V\union W,\gamma_n  \right]\o
A\bigl(
\left(e,V,\sigma \right), \left(u,W,\tau \right)
\bigr)
\bigr>\\ \typeout{Here is C}
=f\o
\bigl<
\Arr \left[e_1,V\union W,\gamma_1  \right]\o
A\bigl(
\left(e_1,V,\gamma_1 \right), \left(u,W,\tau \right)
\bigr),\\ \quad\quad\ldots,
\Arr \left[e_n,V\union W,\gamma_n  \right]\o
A\bigl(
\left(e_n,V,\gamma_n \right), \left(u,W,\tau \right)
\bigr)
\bigr>\end{multline*}
Again a complete proof of the last equality may require a lemma.
Continuing, we have that \typeout{Here is D}
\begin{multline*}\Arr[t]= f\o \bigl< \Arr
\left[e_1[e\leftarrow u], \left(V\setminus\{x\}\right)\union W,
\gamma_1 \right] ,
\ldots,
\Arr \left[e_n[e\leftarrow u], \left(V\setminus\{x\}\right)\union W,
\gamma_n \right]
\bigr>\\ \typeout{Here is E}
\shoveright{\text{(by induction hypothesis)}}\\
=\Arr\bigl[f \left(e_1[x_1\leftarrow u] ,\ldots,
e_n[x_n\leftarrow u] \right),
\left(V\setminus\{x\} \right)\union W,\sigma\bigr]\\ \typeout{Here is F}
\shoveleft{=\Arr\bigl[t\left[x\leftarrow \left(u,W,\tau
\right)\right]\bigr]}\\ \end{multline*}
which is what we get from the
recursive definition.
This completes the proof of the equivalence of the two
definitions.

Later, we shall use the equivalence of these two
methods of obtaining the arrow corresponding to the term in which substitution
has been made. To facilitate reference we record this in the form of a lemma.
\thl\label{facil}
Let $t\ceq \left(e,V,\sigma \right)$ and $t'\ceq \left(u,W,\tau\right)$
be terms, suppose $x\in V$ and suppose $\Outp[u]=\Type[x]$ so
that $u$ may be substituted for $x$. Then there exists an arrow
\begin{equation*}A\ceq A\bigl( \left(e,V,\sigma\right),
\left(u,W,\tau \right)\bigr):
\prd\bigl(\left(V\setminus\{x\}\right)\union W\bigr)
\to \prd\left(V\union W\right)\end{equation*}
so that $$\Arr\bigl[ t
\left[x\leftarrow
u\right],\left(V\setminus\{x\}\right)\union W,\tau \bigr]
= \Arr\bigl[ \left(e,V\union W,\sigma \right)\bigr]\o A$$ \ethl

\section{Rules of inference of MSEL}\label{rulinfms}

In this section we show how the rules of inference of multisorted equational
logic can be codified into our present system.  This is a two-step process.
First, we show that for each rule of inference the pair of equal arrows
corresponding to the conclusion of the rule of inference can be constructed
using the rules of construction of graph-based logic \cite{logstr} from the
single arrow or the product of the equal pairs of arrows that form the
hypothesis of that rule of inference.  Next, we exhibit the construction as an
actual factorization as defined in \cite{logstr}, where the
nodes and arrows appearing in the positions corresponding to the various
labels on the diagram~\eqref{actfact}
are the appropriate instances of the hypothesis,
claim, workspace and so on for the rule in question.
\begin{equation}\label{actfact}\xymatrix{
& \hyp \ar[d]^{\claimcon}\\
\claim \ar[ur]^{\verif}
\ar[r]_{\hypcon}
& \wksp}\end{equation}
While some of these are done in detail some others are not.  For our purposes,
it is enough to prove that a codification as an actual factorization in
$\SynCat\bigl[\FinProd,\Fa\bigr]$ (as defined in
Section~\ref{sksigg}) is possible. In general, this may be
done in more than one way. Symmetry and reflexivity are treated
separately. Transitivity, concretion, abstraction and
substitutivity are all treated in Section~\ref{alltreat}, as
they all are special instances of a worked-out example in
\cite{logstr}.

\subsection{Reflexivity}

The equational rule of inference is
$$\Frac{}{h=_Vh}$$
\astep Translated into the present context, as an instance of the rule
of construction REF \cite{logstr}, this is represented as
$$\Rule{
\xymatrix{ A \ar[r]^f &   B\\
}}{
\xymatrix{ A \rtwo ff & B\\
}}{REF}{}
$$
where
\begin{align*}f&\ceq \Arr\left[h,V,\Type[V]\right]\\
A&\ceq \Dom\bigl[\Arr\left[h,V,\Type[V]\right]\bigr] \\
B&\ceq \Cod\bigl[\Arr\left[h,V,\Type[V]\right]\bigr]\end{align*}
This concludes the first step.

\bstep The corresponding actual factorization:
$$\xymatrix{
&(\arr\x\arr)^{<f,f>}\ar[d]^{\proj_1}\\
\arr^f\ar[r]_{\id}
\ar[ur]^{\Delta}& \arr^f
}$$
Note that one can also use $\proj_2$ as $\claimcon$.
This factorization actually occurs in
$\SynCat[\Cat,\F]$ and is inherited by $\SynCat[\FinProd,\F]$.
A similar remark is true of the constructions for symmetry and
transitivity.

\subsection{Symmetry}

Although in \cite{logstr}, we did not use a rule of construction
corresponding to symmetry, we shall record an actual
factorization for this to facilitate later discussion (in this
section) on proofs as actual factorizations. The rule in
equational deduction is
$$\Frac{e=_Ve'}{e'=_Ve}$$
We define \begin{align*}f&\ceq \Arr\left[e,V,\Type[V]\right]\\
f'&\ceq \Arr\left[e',V,\Type[V]\right]\end{align*}
then the actual factorization is as exhibited below:

$$\xymatrix{
&(\arr\x\arr)^{<g,f>}\ar[d]^{<\Proj[2],\Proj[1]>}\\
(\arr\x\arr)^{<f,g>}
\ar[r]_{\id
\x\id}\ar[ur]^{<\Proj[2],\Proj[1]>}&(\arr\x\arr)^{<f,g>} }$$

\subsection{Transitivity}

The equational rule of inference is
$$\Frac{e=_Ve'\quad e'=_Ve''}{e=_Ve''}$$
For the first step we define
\begin{align*}f:D\rightarrow C&\ceq \Arr\left[e,V,\Type[V]\right]\\
g:D\rightarrow C&\ceq \Arr\left[e',V,\Type[V]\right]\\
h:D\rightarrow C&\ceq \Arr\left[e'',V,\Type[V]\right]\end{align*}
Note that $f$, $g$, and $h$ have the same domain
and the same codomain as $e$, $e'$, and $e''$ have the
same type and as $V$ is the same in each of the terms exhibited
below:

$$\Rule{\xymatrix{ D \ar[r]<1ex>^f \ar[r]<-1ex>_g & C
& D \ar[r]<1ex>^g \ar[r]<-1ex>_h & C\\
}}{\xymatrix{ D \ar[r]<1ex>^f \ar[r]<-1ex>_h & C\\
}}{TRANS}{for all objects $D$ and $C$ and all
arrows $f,g,h:D\to C$ of $\CatTh\bigl[\FinLim,\Sk[\Ssc]\bigr]$}$$

The corresponding actual factorization is provided in
Section~\ref{alltreat}.


\subsection{Concretion}

In this case the equational inference rule reads

Given a set $V$ of typed variables, $x\in V$ and an equation
$e=_Ve'$ such that $x\in
V\setminus(\VarSet[e]\union\VarSet[e'])$, and given that
$\Type[x]$ is inhabited,
$$\Frac{e=_V e'}{e=_{V\backslash \{x\}}e'}$$

We define $\tau\ceq \Type[e]=\Type[e']$ and $\sigma\ceq \Type[x]$,
and
\begin{align*}f:P\rightarrow \tau&\ceq \Arr\left[e,V,\tau\right]\\
f':P\rightarrow \tau&\ceq \Arr\left[e',V,\tau\right]\\
g:Q\rightarrow \tau&\ceq \Arr\left[e,V\backslash\{x\},\tau\right]\\
g':Q\rightarrow \tau&\ceq \Arr\left[e',V\backslash\{x\},\tau\right]\end{align*}

Using Lemmas~\ref{lemmaa} and~\ref{lemmb}, we may choose a map
$$h:\prd\left(V\setminus\{x\}\right)\to\prd V$$ such that
\begin{align*}
g=f\o h\\
g'=f'\o h
\end{align*}
Thus coded as arrows, the rule reads
$$\frac{f=f'}{f\o h=f'\o h}$$

\subsection{Abstraction}

In this case the equational rule of inference reads

Given a set
of typed variables and
$x\in\Vbl[S]\setminus V$, $$\Frac{e=_V
e'}{e=_{V\union\{x\}}e'}$$

We define $\tau\ceq \Type[e]=\Type[e']$ and $\sigma\ceq \Type[x]$,
and
\begin{align*}f:P\rightarrow \tau&\ceq
\Arr\left[e,V,\tau\right]\\ f':P\rightarrow \tau&\ceq
\Arr\left[e',V,\tau\right]\\ g:Q\rightarrow \tau&\ceq
\Arr\left[e,V\cup\{x\},\tau\right]\\ g':Q\rightarrow \tau&\ceq
\Arr\left[e',V\cup\{x\},\tau\right]\end{align*}

Using Lemmas~\ref{lemmaa} and~\ref{lemmb}, we may choose a map
$$h:\prd\left(V\union \{x\}\right)\to V$$ such that
$g=f\o h$ and $g'=f'\o h$.  Thus coded as arrows the rule reads
$$\frac{f=f'}{f\o h=f'\o h}$$

\subsection{Substitutivity}

Given a set
$V$ of typed variables, $x\in V$, and
expressions $u$ and $u'$ for which
$\Type[x]=\Type[u]=\Type[u']$ and $\Type[e]=\Type[e']=\tau$,
$$\Frac{e=_Ve'\quad u=_Wu'}{e[x\leftarrow
u]=_{V\setminus\{x\}\union W}e'[x\leftarrow u']}$$

We already have the representations

\begin{align*}f&\ceq\Arr\left[e,V,\Type[e]\right]=Q[e]I[e]D[e] \\
f'&\ceq\Arr\left[e',V,\Type[e]\right]=Q[e']I[e']D[e']\\
g&\ceq\Arr\left[u,W,\Type[u]\right]=Q[u]I[u]D[u]\\
g'&\ceq\Arr\left[u',W,\Type[u]\right]=Q[u']I[u']D[u']\\
h&\ceq\Arr\bigl[e[x\leftarrow
u],\left(V\setminus\{x\}\right)\union W,\tau\bigr]\\
h'&\ceq\Arr\bigl[e'[x\leftarrow
u'],\left(V\setminus\{x\}\right)\union W,\tau\bigr]
\end{align*}

In view of Lemma~\ref{facil}, we may choose an arrow
$$A:\prd\bigl(\left(V\setminus\{x\}\right)\union
W\bigr)\to\prd\left(V\union W\right)$$ and
$$A':\prd\bigl(\left(V\setminus\{x\}\right)\union
W\bigr)\to\prd\left(V\union W\right)$$
for which $h=f\o A$ and $h'=f'\o A'$.  Note that $A$ and $A'$
are equal, as $f$ and $f'$ and $g$ and $g'$ are (refer to the
definition of $A$ in Lemma~\ref{facil}).

Thus coded in terms of arrows, the rule reads

$$\frac{f=f' \quad\quad A=A'}{f\o A=f'\o A'}$$

\subsection{Transitivity, concretion, abstraction and
substitutivity as  actual factorizations}\label{alltreat}

We recall the following proposition in \cite{logstr}:
\thp
In any category, given the diagram
\begin{equation}\label{twotri}\xymatrix{
A \ar[d]_f \ar[r]^h
& B \ar[d]^k\\
C \ar[ur]|x \ar[r]_g
& D}\end{equation}
if the two triangles commute, then so does the outside square.\ethp
It is shown
there how the proof may be viewed as an actual factorization.

Transitivity may be viewed as a special case of this once equations are
interpreted as commutative diagrams as shown in
Diagram~\eqref{twa}:
\begin{equation}\label{twa}\xymatrix{ D
\ar[d]_{\Id[D]} \ar[r]^h & C \ar[d]^{\Id[C]}\\ D \ar[ur]|g
\ar[r]_f &  C}\end{equation}
The fact that the two triangles commute means that $h=g$ and
$g=f$.  That the outside square commutes means that $h=f$.

In view of Lemmas~\ref{lemmaa} and~\ref{lemmb},
concretion and abstraction can be seen to be special cases of the following:
For every pair of formally equal arrows $f,f':D\to C$ and for
every $h:E\to D$, $f\o h$ and $f'\o h$ are formally equal.  This
can also be realized as a special case of
the commutativity of Diagram~(\ref{twotri}), with choices as
shown: : \begin{equation}\label{twb}\xymatrix{ E \ar[d]_h
\ar[r]^h & D \ar[d]^f\\ D \ar[ur]|{\Id[D]} \ar[r]_g &
C}\end{equation} Particular choices for $h$ yield concretion and
abstraction.

In substitutivity, in view of Lemma~\ref{lemmb}, we have the following in
terms of arrows:  For every pair of formally equal arrows
$f,f':D\to C$, and for every pair of formally equal arrows
$A,A':E\to D$, $f\o A$ and $f'\o A'$ are formally equal.  This
is also a special case of Diagram~(\ref{twotri}) as shown below:
\begin{equation}\label{twc}\xymatrix{
E \ar[d]_{A'} \ar[r]^A
& D \ar[d]^f\\
D \ar[ur]|{\Id[D]} \ar[r]_{f'}
& C}\end{equation}
On the basis of the preceding discussion
we conclude that we may make
choices for all nodes and arrows in the diagram
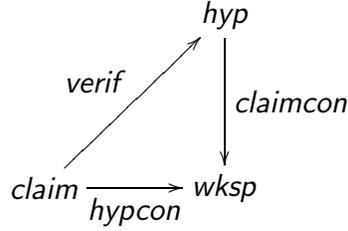
\begin{equation*}\label{actfact2}\xymatrix{ & \hyp
\ar[d]^{\claimcon}\\ \claim \ar[ur]^{\verif} \ar[r]_{\hypcon} &
\wksp}\end{equation*} so that the actual factorization in
$\SynCat[\FinProd,\F]$ codes transitivity, concretion,
abstraction and substitutivity respectively.

\rem Our goal in this section is to produce for every equational
deduction for every equation a corresponding actual
factorization.  However, in order to do this we need to put
equational deductions into some normal form to allow easy
translation.  We shall also need to use certain operations on
actual factorizations.  We record these below as various lemmas.
\erem
\thl\label{evlem}
Every two actual factorizations in any syntactic category
$$\begin{array}{cc}
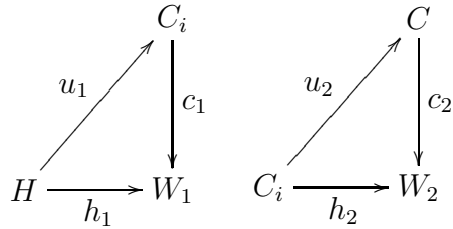
\xymatrix{ & C_i
\ar[d]^{c_1}\\
H \ar[ur]^{u_1} \ar[r]_{h_1}
& W_1} &
\xymatrix{
& C \ar[d]^{c_2}\\
C_i \ar[ur]^{u_2} \ar[r]_{h_2}
& W_2}
\end{array}$$
with the above labels
can be
pasted together to yield a single actual factorization with
labels as shown:
$$
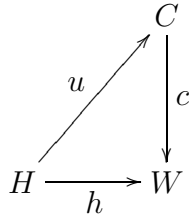
\xymatrix{
& C \ar[d]^{c}\\
H \ar[ur]^{u} \ar[r]_{h}
& W}$$
\ethl
\pf
As every node in a syntactic category (or $\SynCat[\E,\F]$) is the
vertex of a limit cone over some diagram in $\E$,
we may choose
\begin{align*}
\Delta_1 &=\BsDiag[W_1]\\
\Delta_2 &=\BsDiag[W_2]\\
\Delta_i &=\BsDiag[C_i]
\end{align*}
to get the following in the category of diagrams  of
$\E$:

$$ \xymatrix{
& \Delta_1 \ar[d]^{\alpha_1}\\
\Delta_2 \ar[r]_{\alpha_2} & \Delta_i
}$$
where $\alpha_1$ and $\alpha_2$ are the morphisms of
diagrams that give rise to $c_1$ and $h_2$.
As the category
of diagrams in a category is small complete, we may form the
pullback as shown: \begin{equation}\label{pbincd} \xymatrix{
\Delta \ar[r]^{\beta_1} \ar[d]_{\beta_2} & \Delta_1
\ar[d]^{\alpha_1}\\ \Delta_2 \ar[r]_{\alpha_2} & \Delta_i
}\end{equation} Taking the limit over the diagrams corresponding
to the vertices in~\eqref{pbincd} and using the lemmas in
Section~3 of \cite{logstr}, we get the following diagram in
$\SynCat[\E,\F]$:
$$ \xymatrix{
C_i \ar[r]^{h_2} \ar[d]_{c_1}
& W_2 \ar[d]^{d_2}\\
W_1 \ar[r]_{d_1}
& W
}$$
This gives the following diagram in $\SynCat[\E,\F]$:
$$ \xymatrix{
&& C \ar[d]^{c_2}\\
& C_i \ar[ur]^{u_2} \ar[r]_{c_2} \ar[d]^{c_1}
& W_2 \ar[d]^{d_2}\\
H \ar[ur]^{u_1}
\ar[r]_{h_1}
& W_1 \ar[r]_{d_1}
& W}$$
The lemma follows by setting
\begin{align*}
u &\ceq u_2\o u_1\\
c &\ceq d_2\o c_2\\
h &\ceq d_1\o h_1
\end{align*}\epf
\rem The above is an analogue of getting the deduction
$\frac{E_1}{E}$ given the deduction $\frac{E_1}{E_2}$ and
$\frac{E_2}{E}$.\erem

\thl\label{givlem}
Given two actual factorizations in any syntactic category
$$\begin{array}{cc}
\xymatrix{ & C_1
\ar[d]^{c_1}\\
H_1 \ar[ur]^{u_1} \ar[r]_{h_1}
& W_1} &
\xymatrix{
& C_2 \ar[d]^{c_2}\\
H_2 \ar[ur]^{u_2} \ar[r]_{h_2}
& W_2}
\end{array}$$
we have the factorization
$$ \xymatrix{
& C_1\x C_2 \ar[d]^{c_1\x c_2}\\
H_1\x H_2 \ar[ur]^{u_1\x u_2} \ar[r]_{h_1\x h_2}
& W_1\x W_2
}$$\ethl
\pf Omitted.\epf

\subsection{Normal forms for equational deductions}

We recall the usual recursive definition of a deduction of the
equation $E$ from the family of equations $(E_i\mid i\in
1\twodots n)$.  A deduction is a tree with $E$ at the root and
$n$ nodes
$(E_i\mid i\in
1\twodots n)$
at level $0$.  For all $m\in 1\twodots (n-1)$, the nodes are all
obtained by the rules of inference (listed in
Section ~\ref{rulinfms}) from either (a) any one or two nodes at
strictly smaller levels (as the rules of inference that have a
nonempty set for a premise have either one or two premises), or
(b) the empty premise.

We assume given a deduction $$D=\bigl(\left(E_{lw}\mid w\in
0\twodots W(l)\right)\mid i\in 0\twodots L\bigr)$$
where $l$ refers to the level in the deduction tree and $W(l)$
refers to the width of the deduction tree at level $l$.  In
addition: (a) $W(0)=n$, that is, the list of premises has $n$
entries $\left(E_i\mid i\in 1\twodots n\right)$, and (b)
$W(L)=1$ and has exactly one entry, namely $E$.

In order to construct an actual factorization corresponding to
$D$, we first transform $D$ into a normal form to
facilitate coding. The normal form described below will involve
considerable redundancy.  We list below the relevant features of
the normal form. \blist{NF} \item For all $l\in 1\twodots L$ and
for all $w\in 0\twodots W(l)$, the premises from which $E_{lw}$
are deduced by one of the rules of inference
in Section ~\ref{rulinfms} appear at the
immediately preceding level, that is level $l-1$.  (Note that
$l\in 1\twodots L$, in other words we exclude the hypotheses at
level $0$.) \item For all $l\in 0\twodots L$ and for all $w\in
0\twodots W(l)$, $E_{lw}$ is used exactly once for a deduction
of some $E_{l+1,w}$ at the immediately following level. \elist
NF.1 is achieved by carrying over every single $E_{lw}$ at level
$l$ to the root level, using the rule $\frac{E}{E}$. NF.2 is
achieved by repeating every hypothesis starting at level $0$
with every equation in the list $\left(E_i\mid i\in 1\twodots
n\right)$ as many times as it is used in the body of the given
deduction $$\bigl(\left(E_{lw}\mid w\in 0\twodots
W(l)\right)\mid i\in 0\twodots L\bigr)$$ We shall give each such
deduction the name $$\Ddcn\bigl[ \left(E_i\mid i\in 1\twodots
n\right), E\bigr]= \bigl(\left(E_{iw}\mid w\in 0\twodots
W(l)\right)\mid l\in 0\twodots L\bigr)$$ We record the preceding
discussion in the following Lemma. \thl\label{lastlem} Given any
deduction of $E$ from the family of equations $\left(E'_i\mid
i\in 1\twodots n'\right)$, we may find an equivalent one (in
normal form) $$\Ddcn\bigl[\left(E_i\mid i\in 1\twodots
n\right),E\bigr]= \bigl(\left(E_{ln}\mid w\in 0\twodots
W(l)\right)\mid l\in 0\twodots L\bigr)$$ that has the following
properties: \begin{enumerate}

\item For all $l\in 1\twodots L$ and for all $w\in 0\twodots
W(l)$, $E_{lw}$ follows from some members of the list
$\left(E_{l-1,w}\mid w\in 1\twodots W(l-1)\right)$ or the empty
premise using some rule of inference.
\item For all $l\in 0\twodots L$ and for all $w\in 0\twodots
W(L)$, $E_{lw}$ is used exactly once for deducing some
$E_{l+1,w}$ in the next level using some rule of inference.

\end{enumerate}

\ethl

\thm
Given any deduction $\Ddcn\bigl[\left(E_i\mid e\in 1\twodots
n\right),E\bigr]$ of equational logic, we may construct an
actual factorization in $\SynCat[\FinProd,\F]$ that
corresponds to it.
\ethm

\pf
In view of Lemma~\ref{lastlem}, we may, without loss of
generality, assume that the given deduction is in normal form,
that is,
$$\Ddcn\bigl[\left(E_i\mid i\in 1\twodots n, E\right), E\bigr]
=\bigl(\left(E_{lw}\mid w\in 0\twodots W(l)\right)\mid l\in
0\twodots L\bigr) $$
Because the deduction is in normal form for every pair of
consecutive levels $(l, l+1)\in(0\twodots l)\x(1\twodots l)$ we
may choose partitions of the index sets $0\twodots W(l)$ and
$0\twodots W(l+1)$ with the following properties:
\begin{enumerate}

\item Both index sets have the same number $P(l,l+1)$ of parts.
For all $p\in 0\twodots P(l,l+1)$, we define $p_l$ to be the
family of equations in the $p$th part at level $l$ and similarly
$p_{l+1}$ to be the family of equations in the $p$th part at
level $l+1$.

\item For all $p\in 0\twodots P(l,l+1)$, $\frac{p_l}{p_{l+1}}$ is
an instance of some rule of inference; that is, the pair of
levels $(l, l+1)$ may be rewritten as
$$\bigl(\frac{p_l}{p_{l+1}}\mid p\in 0\twodots
P\left(l,l+1\right)\bigr)$$
\end{enumerate}

As we have already shown how each rule of inference can be coded
as an actual factorization in Section~\ref{alltreat} it follows
that we may select, for each pair of levels $(l,l+1)\in(0\twodots
L)\x(1\twodots L)$,
a family (indexed by $p\in 0\twodots P(l,l+1)$)
of actual factorizations of the form
\begin{equation}\label{famm}\xymatrix{
& C_{p,l} \ar[d]^{c_{p,l}}\\
H_{p,l} \ar[ur]^{u_{p,l}} \ar[r]_{h_{p,l}}
& W_{p,l}
}\end{equation}
in $\SynCat[\FinProd,\F]$,
where each instance of~\eqref{famm}
corresponds to $$\frac{p_l}{p_{l+1}}$$ via the coding of
inference rules.

Using Lemma~\ref{givlem}, we may combine for all $l\in
0\twodots(L-1)$ this family of actual factorizations into a
single actual factorization
$$ \xymatrix{
& C_l \ar[d]^{c_{p,l}}\\
H_l \ar[ur]^{u_l} \ar[r]_{h_{p,l}}
& W_l
}$$
where for all $x\in\{H,C,W,u,c,h\}$, $$x_l=\prd_{p\in P(l,l+1)}
x_{p,l}$$
Note that $P(l-1,l)$ and $P(l,l+1)$ may be different.  However,
for all $l\in 0\twodots(l-1)$, $C_l$ and $H_{l+1}$
will be isomorphic because of associativity.  For every $l$ we
define $\alpha_{l,l+1}$ to be the associativity isomorphism.  We
now have the following diagram.
$$
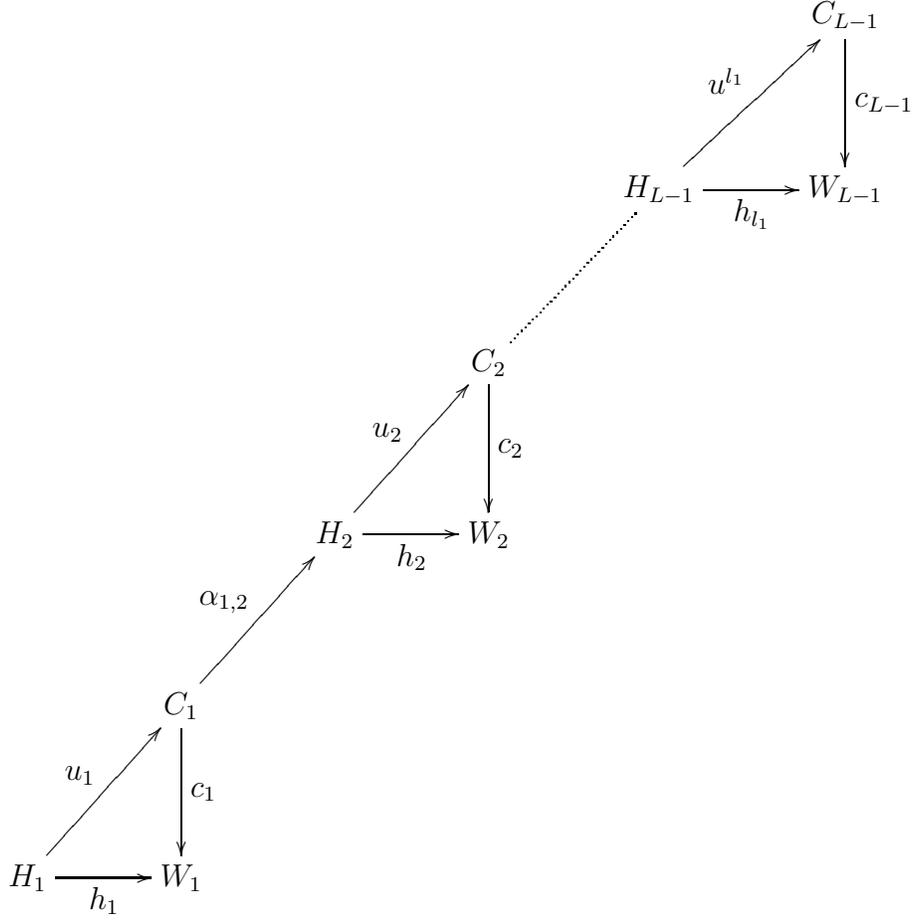
\xymatrix{
&&&&& C_{L-1} \ar[d]^{c_{L-1}}\\
&&&& H_{L-1} \ar[ur]^{u^{l_1}} \ar[r]_{h_{l_1}}
& W_{L-1}\\
&&&C_2 \ar[d]^{c_2} \ar@{.}[ur]\\
&& H_2 \ar[ur]^{u_2} \ar[r]_{h_2}
& W_2\\
& C_1 \ar[ur]^{\alpha_{1,2}} \ar[d]^{c_1}\\
H_1 \ar[ur]^{u_1}\ar[r]_{h_1}
& W_1
}
$$
Using Lemma~\ref{evlem}, we can combine all of this into a
single actual factorization
\begin{equation}\label{A}
\xymatrix{
& C_{L-1} \ar[d]^{c'}\\
H_1 \ar[ur]^{u'} \ar[r]_{h'}
& W_{L-1}
}
\end{equation}

Finally we note that we have the commutative diagram (in
$\CatTh[\FinProd,\S]$)
\begin{equation}\label{B}
\xymatrix{ C_{L-1} \ar[r]^{\proj}
\ar[d]_{c_{L-1}} & C \ar[d]^{c_{C,W}}\\ W_{L-1} \ar[r]_{\proj} &
W } \end{equation} where
$$ \xymatrix{ & C \ar[d]_{c_{C,W}}\\ H
\ar[r]^{h_{C,W}} &W } $$
is the actual factorization
corresponding to the deduction of $E$ from the corresponding
premises in the partition in level $L-1$.  We note that $E$ is
the only member of one of the parts of level $L$.  Putting
\eqref{A} and~\eqref{B} together, we have
\begin{equation}\label{C}
\xymatrix{& C \ar[d]^c\\
H_1\ar[ur]^u \ar[r]_L
&W}
\end{equation}
where
\begin{align}
u &= \proj\o u'\\
h &= \proj\o h'\\
c&=c_{C,W}
\end{align}
Diagram~\eqref{C} is the actual factorization corresponding to
the given deduction. \epf

\rem
The above shows that every deduction in MSEL occurs as an actual
factorization in $\SynCat[\FinProd,\F]$. \erem



\rem
Although we worked out the details for multisorted equational
logic  and $\CatTh[\FinProd,\F]$, the method will work for any
logical system that can be described as a constructor-space
sketch.  Thus, in general, we shall have some logical system $L$
and a category $\CatTh[\E_{L},\F]$  in which $\E_L$ is the kind
of category in which the models of $L$ are.  For instance, if
$L$ is the typed $\lambda$-calculus, $E_L$ would be $\CCC$, and
if $L$ is intuitionistic type theory, then $E_L$ would be a
constructor space for toposes.

Given any sound and complete deductive system for $L$, if we
interpret terms as arrows and encode them in $\CatTh[\E_L,\F]$
as we have done here, then we conjecture that the method will
show that all theorems of $L$ can be realized as actual
factorizations in $\CatTh[\E_L,\F]$.  (Indeed it appears nearly
obvious that this will happen if we know that $L$ and $E_L$ have
equivalent models; a detailed proof, is of course necessary to
clinch the matter.) In the examples of the preceding paragraphs,
we might use the deductive systems formulated in
\cite{lambekscottbook}. The method used here is quite general.
\erem

\section{Acknowledgments}
The diagrams were prepared
using K. Rose's {\tt xypic}.

\tolerance=4000
\raggedright

\nocite{powerwells}
\nocite{gensk}
\nocite{lambekscottbook}

\end{document}